\documentclass{amsart}
\usepackage{amsmath}
\usepackage{amsfonts}
\usepackage{amssymb}
\usepackage{latexsym}
\usepackage{enumerate}
\newcommand{\restrict}{\upharpoonright}
\newcommand{\eop}{\Box}

\def\cl{c\ell}
    \newtheorem{theorem}{Theorem}[section]
  \newtheorem{lemma}[theorem]{Lemma}
  \newtheorem{corollary}[theorem]{Corollary}
  \newtheorem{proposition}[theorem]{Proposition}

 \newtheorem{*theorem}[theorem]{*Theorem}

  \newcounter{rmk}[section]
  \renewcommand{\thermk}{\thesection.\arabic{rmk}}
  \newenvironment{remark}{\refstepcounter{rmk}
    \medskip\par\noindent {\bf Remark \thermk \enspace}}{
    \par\medskip}

  \newcounter{defn}[section]
 \renewcommand{\thedefn}{\thesection.\arabic{defn}}
  \newenvironment{definition}{\refstepcounter{defn}\bigskip
    \noindent{\bf Definition \thedefn \enspace
      }}{\par\bigskip}

  \newcounter{xmpl}[section]
  \renewcommand{\thexmpl}{\thesection.\arabic{xmpl}}
  \newenvironment{example}{\refstepcounter{xmpl}
    \medskip\par\noindent {\bf Example \thexmpl
      \enspace}}{ \par\medskip}

 \newenvironment{notation}{\bigskip \noindent{\bf Notation.\
}}{\par\bigskip}

 \newenvironment{claim}{\medskip \par\noindent{\sl Claim.\
}}{\par\medskip}

\def\aa{\bar a}
\def\bb{\bar b}
\def\cc{\bar c}
\def\dd{\bar d}
\def\ee{\bar e}
\def\ff{\bar f}

\def\vv{\bar v}

\def\xx{\bar x}
\def\yy{\bar y}
\def\zz{\bar z}
   \def\dnfo{\;\raise.2em\hbox{$\mathrel|\kern-.9em\lower.4em\hbox
{$\smile$}$}}

\def\dnf#1{\lower.9em\hbox{$\buildrel\dnfo\over{ \scriptstyle  #1}$}}

\def\dfo{\;\raise.2em\hbox{$\mathrel|\kern-.9em\lower.4em\hbox{$\smile$}
\kern-.72em\lower.07em\hbox{\char'57}$}\;}
\def\df#1{\lower1em\hbox{$\buildrel\dfo\over{\scriptstyle #1}$}}

\def\indep#1#2#3{\hbox{\mathsurround=0pt$#1 \  \dnf{#2}  \ #3$}}

\newcommand{\cR}{\mathcal{R}}
\newcommand{\Saut}[1]{\operatorname{Saut}_{#1}(M)}
\newcommand{\aut}[1]{\operatorname{Aut}_{#1}(M)}

\begin{document}
\title{Simple homogeneous models}
\author{Steven Buechler}
\thanks{Research of the first author partially supported by the NSF}
\address{Steven Buechler,
Department of Mathematics,
255 Hurley Hall,
University of Notre Dame,
Notre Dame, IN 46556}
\author{Olivier Lessmann}
\address{Olivier Lessmann,
Mathematical Institute,
24-29 St. Giles,
 Oxford University,
 Oxford OX1 3LB,
 United Kingdom}
\subjclass[2000]{03C45}
\keywords{stability theory, simple theories, nonelementary classes}
\date{July 25, 2002}
\maketitle

\begin{abstract}
Geometrical stability theory is a powerful set of model-theoretic
tools that can lead to structural results on models of a simple
first-order theory. Typical results offer a characterization of
the groups definable in a model of the theory. The work is carried
out in a universal domain of the theory (a saturated model) in
which the Stone space topology on ultrafilters of definable
relations is compact. Here we operate in the more general setting
of  homogeneous models, which typically have noncompact Stone
topologies. A structure $M$ equipped with a class of finitary
relations $\cR$ is \emph{strongly $\lambda-$homogeneous} if orbits
under automorphisms of $(M,\cR)$ have finite character in the
following sense: Given $\alpha$ an ordinal $<\lambda\leq|M|$ and
sequences $\aa=\{\,a_i:\:i<\alpha\,\}$,
$\bb=\{\,b_i:\:i<\alpha\,\}$ from $M$,  if
$(a_{i_1},\dots,a_{i_n})$ and $(b_{i_1},\dots,b_{i_n})$ have the
same orbit, for all $n$ and $i_1<\dots<i_n<\alpha$, then
$f(\aa)=\bb$ for some automorphism $f$ of $(M,\cR)$. In this paper
strongly $\lambda-$homogeneous models $(M,\cR)$ in which the
elements of $\cR$ induce a symmetric and transitive notion of
independence with bounded character are studied. This notion of
independence, defined using a combinatorial condition called
``dividing'', agrees with forking independence when $(M,\cR)$ is
saturated. The concept central to the development of geometrical
stability theory for saturated structures, namely the canonical
base, is also shown to exist in this setting. These results
broaden the scope of the methods of geometrical stability theory.
\end{abstract}

This paper attempts to give a self-contained development of
dividing theory (also called forking theory) in a strongly
homogeneous structure. Dividing is a combinatorial property on the
invariant relations on a structure that have yielded deep results
for the models of so-called ``simple'' first-order theories. Below
we describe for the nonspecialist how this paper fits in the
broader context  of geometrical stability theory. Naturally, some
background in first-order model theory helps to understand these
motivating results, however virtually no knowledge of logic is
assumed in this paper. Readers desiring a more thorough
description of geometrical stability theory are referred to the
surveys \cite{hr:current} and \cite{hr:icm}.

Traditionally, geometrical stability theory is a collection of
results that apply to definable relations on arbitrary models of
a complete first-order theory. It is equivalent and convenient to
restrict our attention to the definable relations on a fixed
representative model of the theory, called a universal domain.
Using the terminology of the abstract, a \emph{universal domain}
is an uncountable model $M$ equipped with the first-order
definable relations $\cR$ which is strongly $|M|-$homogeneous and
\emph{compact}; i.e., if $\{\,X_i:\:i\in I\,\}$, where $|I|<|M|$,
is a family of definable relations on $M$ so that $\bigcap_{i\in
F}X_i\not=\emptyset$ for any finite $F\subset I$, then
$\bigcap_{i\in I}X_i\not=\emptyset$. For our purposes the reader
can assume there is a one-to-one correspondence between (complete
first-order) theories and universal domains.

A massive amount of abstract model theory was developed en route
to Shelah's proof of Morley's Conjecture about the number of
models, ranging over uncountable cardinals,
 of a fixed first-order theory \cite{sh:book}.
Most of the work concerned the case of a stable theory, which will
not be defined here for the sake of brevity. What is relevant is
that most theorems describing the models of a stable theory rely
on the forking independence relation. The forking independence
relation, $F$, is a ternary relation on the subsets of the
universal domain of a theory, where $F(A,B;C)$ is read ``$A$ is
forking independent from $B$ over $C$'' (see
Remark~\ref{rem:defn_free}). In a stable theory forking
independence is symmetric (over $C$), has finite character (in $A$
and $B$), bounded dividing, the free extension property and is
transitive. (See Definition~\ref{defn:simple},
Theorem~\ref{thm:symm_lemma} and Corollary~\ref{cor:trans} for
precise statements of these properties.) These properties
facilitate the introduction of several notions of dimension that
lead to procedures for determining when two models are isomorphic.
The combinatorial-geometric properties of the definable relations
reflected in these dimensions profoundly impact the structure of
the models beyond the question of fixing an isomorphism type. The
results connected to algebra, known as geometrical stability
theory, have lead to new theorems in number theory, in particular
the first proof of the Mordell-Lang Conjecture for function fields
in positive characteristics.

In hind sight the origins of geometrical stability theory can be
traced to the work of Zil'ber \cite{zil:unc}, Chelin, Harrington
and Lachlan \cite{chl} in the late 1970s and early 1980s. The
impact of the area escalated in the mid-1980s with Hrushovski's
discovery that natural algebraic objects (groups, fields, vector
spaces, etc.) are definable  in a universal domain satisfying
various abstract model-theoretic hypotheses. This interplay
between the combinatorial geometry of abstract stability theory
and the model theory of algebraic objects significantly deepened
our understanding of stable theories and gave new insights into
differential fields \cite{pi:diff_galI}, \cite{pi:diff_galII},
\cite{ma.pi:diff_galIII} and even algebraically closed fields
\cite{hr.zi:zariski}.

Looking towards applications to number theory, Hrushovski realized
the need for an analysis of a ``generic'' algebraically closed
difference field (i.e., an algebraically closed field adjoined
with an automorphism that is in some sense universal for
difference fields). While the relevant universal domain is clearly
unstable, Chatzidakis and Hrushovski showed \cite{ch.hr:diff} that
it does admit a notion of independence reminiscent of forking.
After this work was well underway it was realized that their
notion of dependence agrees with forking independence and the
relevant universal domain is ``simple''. In \cite{sh:simp} Shelah
defines a simple theory in terms of a combinatorial property on
families of definable sets and proves that forking independence
has some of the nice properties found in a stable theory. Every
stable theory is simple but not conversely. Our understanding of
simple theories increased dramatically when Byunghan Kim showed in
his doctoral research \cite{ki:thesis} that forking in a simple
theory is symmetric, transitive, and satisfies type amalgamation
(which compensates for the loss of the definability of types, true
in a stable theory but not in a simple unstable theory). Since
Kim's seminal work more of the machinery of geometrical stability,
such as canonical bases and generics in groups and fields, has
been generalized to simple theories (\cite{hakipi},
\cite{bu.pi.wa:supersimp}, \cite{wa:hdef_gps}). It is hoped that
the geometrical stability theory of simple theories will have
additional applications to number theory.
\bigskip

All of the model theory discussed above takes place in the context
of a universal domain of a first-order theory. Its applicability
to number theory depends on the first-order axiomatizability of
fundamental concepts in algebraic geometry. In this paper we begin
laying the foundation for the application of geometrical stability
theory to some mathematical areas that cannot be captured with
first-order logic. Whereas Kim's development of forking takes
place in a universal domain, the context here is a strongly
homogeneous model. That is, we drop the requirement that the model
is compact. Using a definition of simplicity for a strongly
homogeneous model that specializes to simplicity for a universal
domain, it is shown here that forking independence satisfies the
same basic properties (symmetry and transitivity, e.g.) as in a
simple universal domain. Moreover, the higher-order theorems like
type amalgamation and parallelism of types that are critical to
geometrical stability theory also hold in the simple strongly
homogeneous setting. In Section~\ref{sec:examples} examples are
given of strongly homogeneous models that are simple, although
they are not models of simple theories. In particular, it is shown
that any Hilbert space is a subspace of a strongly homogeneous
Hilbert space and the latter is simple. In \cite{be.bu:hilbert}
Berenstein lays the groundwork for deeper applications of model
theory to functional analysis by showing that many structures of
the form $(H,T)$, where $H$ is a strongly homogeneous Hilbert
space and $T$ is a bounded linear operator on $H$, are strongly
homogeneous and simple. Forking independence in a Hilbert space is
equivalent to a notion of independence based on orthogonal
projection, a very natural geometrical relation. At this early
``proof of concept'' stage it is difficult to gauge the potential
impact of this work on the understanding of these examples from
analysis. Berenstein also shows in \cite{be:thesis} that analogues
of theorems from stable group theory extend to the homogeneous
setting.

This paper is far from the first investigation of
stability-theoretic concepts in models that are not universal
domains. Shelah, Grossberg, Hyttinen, Lessmann and others have
extensively studied the spectrum function of a class of models
that is not the class of models of first-order theory. (The
spectrum function of a class models assigns to a cardinal number
$\lambda$ the number of models in the class of cardinality
$\lambda$, up to isomorphism.) Analogues of forking independence
are found in many of these papers, especially \cite{sh.hy:strong},
although fewer properties can be expected in this setting than for
forking independence in a simple first-order theory.
>From a hypothesis about the spectrum function of a class of models it is
normally impossible to show that the class contains a strongly
homogeneous model. Thus, from the perspective of that body of
research, the context of this paper is very limiting. However,
there are very natural mathematical objects that are strongly
homogeneous and simple although not the models of simple
first-order theories, making this work worthwhile.

This work was strongly influenced by \cite{pi:exist} and its
precursor \cite{hr:robinson}. The proofs of the basic properties
of forking independence in those papers showed that Kim's
treatment could be reproduced in a non first-order setting.

\section{Logical structures and homogeneity}

\subsection{Definition of a logical structure}

Let $M$ be a model in a language $L$. That is $M$ is a set
together with  finitary relations and functions corresponding to
symbols in $L$.

\begin{definition}
\label{defn:logical} The pair $(M,\cR)$ is a \emph{logical
structure} if $M$ is a structure for a first-order language $L$
and $\cR$ is a collection of finitary relations on $M$ satisfying
the following requirements.
\begin{enumerate}[(1)]
\item The elements of $\cR$ are invariant under automorphisms of
$M$. \item $\cR$ is closed under finite unions and intersections.
\item If $R\in\cR$ is $n-$ary and $\pi$ is a permutation of $n$,
then $\{\,(a_{\pi(1)},\dots,a_{\pi(n)}):\:R(a_1,\dots,a_n)\,\}$,
is also in $\cR$.
\end{enumerate}
\end{definition}

\begin{remark}
For the reader unfamiliar with first-order languages an equivalent
formulation can be used. Instead of $M$ being a structure in
language we consider a faithful group action $(G,M)$ and view $G$
as the automorphism group of the structure. Then $\cR$ is a class
of $G-$invariant relations satisfying (2) and (3).
\end{remark}

\begin{definition}
\label{defn:rtp} Let $\aa\in A^\alpha$, $\bb\in A^\beta$, where
$\alpha$, $\beta$ are ordinals. For $i<\alpha$, let $v_i$ be a
variable that ranges over the elements of $M$. Let
$tp_{\cR}(\aa/\bb)$, called the \emph{$\cR-$type of $\aa$ over
$\bb$ in $M$}, be the set of expressions
$R(v_{i_1},\dots,v_{i_m},b_{j_1},\dots,b_{j_l})$, where $R$ is an
$m+l$-ary relation in $\cR$ and
$R(a_{i_1},\dots,a_{i_m},b_{j_1},\dots,b_{j_l})$. If
$\bb=\emptyset$ it may be omitted.

A set $p$  is an \emph{$\cR-$type in $\vv$ over $A$} if it
consists of expressions of the form
$R(x_1,\dots,x_n,a^R_1,\dots,a^R_m)$, where $R\in\cR$,
$x_1,\dots,x_n$ are variables from $\vv$ and $a^R_1,\dots,a^R_m\in
A\subset M$. An $\cR-$type is \emph{consistent} if it is realized
in $M$. If $p(\vv)$ is an $\cR-$type in $\vv$ over $\bb$, the
sequence $\cc$ \emph{realizes} $p$ if $tp_{\cR}(\cc/\bb)\supset
p$, and $p(M)$ denotes $\{\,\cc:\:\cc\textrm{ realizes }p\,\}$. An
$\cR-$type $p$ over $A$ is \emph{complete} if for some $\aa$ and
$q=tp(\aa/A)$, $p(M)=q(M)$.

If the class of relations $\cR$ is clear from context we will drop
it from the notation and  write $M$ for $(M,\cR)$ and $tp(a/b)$
for $tp_{\cR}(a/b)$.
\end{definition}

Looking ahead, dividing will be defined for consistent $\cR-$types
over subsets of $M$.

\begin{example}
\label{xmpl:logic_struct} (i) Let $K$ be a field and $\cR$ the
collection of all constructible sets on $K$. Then $(M,\cR)$ is a
logical structure.

(ii) Let $M$ be a structure in the language $L$ and $\cR$ the
class of first-order definable relations on $M$. Then, $(M,\cR)$
is a logical structure.

(iii) Let $M$ be a structure in the language $L$. Then $(M,\cR)$
is a logical structure if $\cR$ is the collection of all sets of
realizations of types in one of the following families:
\begin{itemize}
\item formulas of first-order logic; \item quantifier-free
formulas of first-order logic; \item finite unions and
intersections of complete first-order types, quantifier-free types
or existential types; \item types in the logic $L_{\kappa\omega}$,
where $\kappa$ is infinite; \item types in $L_k$, the logic with
only $k<\omega$ variables.
\end{itemize}

(iv) If $R$ is a real-closed field and $\cR$ is the collection of
semi-algebraic sets over $R$, then $(R,\cR)$ is a logical
structure.
\end{example}

\begin{remark}
The definition of a logical structure is designed to encompass all
of the examples in (iii). This prohibits us from requiring $\cR$
to be closed under projection, composition or negation.
\end{remark}

\subsection{Homogeneous logical structures}

In this paper our attention will focus on the following kind of
objects.

\begin{definition}
\label{defn:homog} Let $(M,\cR)$ be a logical structure, where $M$
is infinite and $\lambda\leq |M|$ is infinite. Then $(M,\cR)$ is
\emph{strongly $\lambda-$homogeneous} if for all $\alpha<\lambda$
and $\aa,\,\bb\in A^\alpha$, if $tp(\aa)=tp(\bb)$, then there is
an automorphism $f$ of $(M,\cR)$ with $f(\aa)=\bb$.

For brevity we write \emph{$s\lambda-$homogeneous} for strongly
$\lambda-$homogeneous.
\end{definition}

When $\cR$ is the collection of first-order definable relations,
the $s\lambda-$homogeneity of $(M,\cR)$ is equivalent to $M$ being
$s\lambda-$homogeneous as a first-order structure. The following
is trivial but helps to connect $s\lambda-$homogeneity to a more
familiar concept.

\begin{proposition}
\label{prop:homog_auto} An $s\lambda-$homogeneous logical
structure $(M,\cR)$ is $\lambda-$homogeneous in the sense that for
$\alpha< \lambda$,
 $\aa,\,\bb\in A^\alpha$ with
$tp(\aa)=tp(\bb)$ and $c\in M$, there is $d\in M$ such that
$tp(\aa c)=tp(\bb d)$.
\end{proposition}

Strongly $\lambda-$homogeneous models are ubiquitous in
first-order model theory:

\begin{lemma}
\label{lem:s-homog_exist} For $\lambda$ an infinite cardinal and
$T$ a complete first-order theory of cardinality $\leq\lambda$,
there is an $s\lambda-$homogeneous model $M$ of $T$ of cardinality
$\leq 2^\lambda$.
\end{lemma}

\noindent This is \cite[Proposition 2.2.7]{bu:book}.

Homogeneity can be used to obtain the consistency of unions of
chains of complete types. Such a result is most commonly proved
with compactness, however, it also holds in this setting.

\begin{lemma}
\label{lem:union_chain} Let $(M,\cR)$ be an $s\lambda-$homogeneous
logical structure, $A\subset M$, $|A|<\lambda$. Let $A_0\subset
A_1\subset\dots\subset A_\alpha\subset\dots$, $\alpha<\beta$,
 be a chain of sets with $A=\bigcup_{\alpha<\beta}A_\alpha$ and
$p_\alpha(\vv)\in S(A_\alpha)$, $|\vv|<\lambda$, such that
$p_\alpha\subset p_\gamma$, for any $\alpha<\gamma<\beta$. If each
$p_\alpha$, $\alpha<\beta$ is consistent in $M$, then
$\bigcup_{\alpha<\beta}p_\alpha$ is consistent in $M$.
\end{lemma}

\begin{proof}
This follows quickly from \cite[Lemma 5.1.18]{ck:book}.
\end{proof}

\begin{definition}
\label{defn:large} Given an  $s\lambda-$homogeneous logical
structure $(M,\cR)$, an $\cR-$type $p$ is called \emph{large} if
the set of realizations of $p$ has cardinality $\geq\lambda$.
\end{definition}

\begin{definition}
\label{defn:indisc} Let $(M,\cR)$ be a logical structure, $X$ an
ordered set, $A\subset M$ and $I=\{\,a_i:\:i\in X\,\}$ a set of
sequences from $M$ indexed by $X$. Then, $I$ is called
\emph{$A-$indiscernible} or \emph{indiscernible over $A$} if for
all $n<\omega$, and $i_1<\dots<i_n$ and $j_1<\dots<j_n$ from $X$,
$tp(a_{i_1},\dots,a_{i_n}/A) =tp(a_{j_1},\dots,a_{j_n}/A)$.

If $I$ is $A-$indiscernible the \emph{type diagram} of $I$ over
$A$ is the collection of all types $tp(\aa/A)$, as $\aa$ ranges
over finite sequences from $I$ whose indices are increasing in
$X$.
\end{definition}

\begin{remark}
\label{rem:defn_indisc} Let $(M,\cR)$ be an
$s\lambda-$homogeneous logical structure, $X$ an ordered set,
$A\subset M$ and $I=\{\,a_i:\:i\in X\,\}$, $J=\{\,b_i:\:i\in
X\,\}$ $A-$indiscernible sequences with the same type diagram over
$A$ such that $|\cup I|<\lambda$. Then, there is an automorphism
$f$ of $(M,\cR)$ fixing $A$ such that $f(a_i)=b_i$, for $i\in X$.
\end{remark}

The existence and ubiquity of indiscernible sequences in
$s\lambda-$homogeneous models will play a big role in this study.
The following result guarantees that when $\lambda$ is
sufficiently large the set of realizations of a large type over a
relatively small set will contain an infinite sequence of
indiscernibles.

\begin{lemma}
\label{lem:exist_indisc} Let $(M,\cR)$ be $s\lambda-$homogeneous
and $\lambda_0$ be the cardinality of the set of complete
$\cR-$types over $\emptyset$ in finitely many variables realized
in $M$. For each $\lambda_1$ and $\lambda$  sufficiently large,
there is a $\lambda_2\leq\lambda$ (depending on $\lambda_1$ and
$\lambda_0$) such that if $A\subset M$, $|A|\leq\lambda_1$, $X$ is
an ordered set of size at least $\lambda_2$  and
$I=\{\,\aa_i:\:i\in X\,\}\subset M$, where $|\aa_i|\leq\lambda_1$,
for $i\in X$, \emph{then} there is
$J=\{\,\bb_i:\:i<\omega\,\}\subset M$, indiscernible over $A$ such
that for every $n<\omega$, there exists $i_0<\dots<i_n$ in $X$
with
\[
tp(\bb_0,\dots,\bb_n/A)=tp(\aa_{i_0},\dots,\aa_{i_n}/A).\]
\end{lemma}

\begin{proof}
Considering $M$ as an ordinary first-order structure in a language
$L_0$ in which it has elimination of quantifiers, let $T_1$ be an
expansion of $Th(M)$ with Skolem functions, $L_1$ the language of
$T_1$. Let $\lambda_2=\beth_{(2^{\lambda_1+\lambda_0})^+}$. Using
a standard application of the Erdos-Rado Theorem there exists in
some model $N$ of $T_1$, a sequence $\{\,\dd_i:\:i<\omega\,\}$
indiscernible over $A$ in $T_1$ such that for every $n<\omega$,
there exists $i_0<\dots<i_n$ in $X$ with
\[
tp_{L_1}(\dd_0,\dots,\dd_n/A)=tp_{L_1}(\aa_{i_0},\dots,\aa_{i_n}/A).\]
Without loss of generality, $N$ is the Skolem hull of $A\cup J$,
hence every complete type in $L_0$ in finitely many variables
realized in $N$ is also realized in $M$. Assuming that
$|N|<\lambda$, by the $s\lambda-$homogeneity of $M$ there is
$\{\,\bb_i:\:i<\omega\,\}\subset M$, such that
\[
tp(\bb_0,\dots,\bb_n/A)=tp(\dd_0,\dots,\dd_n/A).\] This completes
the proof.
\end{proof}

As a consequence of the preceding lemma, when $M$
$s\lambda-$homogeneous and $\lambda$ is sufficiently large, sets
of realizations of large types over small sets contain infinite
indiscernible sequences. The converse follows from the next lemma.

\begin{lemma}
\label{lem:large_indisc} Let $(M,\cR)$ be $s\lambda-$homogeneous,
$A\subset M$ with $|A|<\lambda$, and suppose $I=\{\,\aa_i:\:i\in
X\,\}\subset M$ is an indiscernible sequence over $A$, where
$|\aa_0|<\lambda$ and $|X|<\lambda$.

(i) For any ordered set $Y$ extending $X$ with $|Y|\leq\lambda$,
there are sequences $\aa_j\in M$, for $j\in Y$, such that
$J=\{\,\aa_i:\:i\in Y\,\}$ is $A-$indiscernible. A fortiori,
$tp(\aa_0/A)$ is large.

(ii) For any linearly ordered set $X'$ of cardinality
$\leq\lambda$ there is $J=\{\,\bb_i:\:i\in X'\,\}\subset M$
indiscernible over $A$ with the same diagram over $A$ as $I$.
\end{lemma}

\begin{proof}
(i)  follows from Lemma~\ref{lem:union_chain} and (ii) follows
from (i) by taking $Y=X+X'$.
\end{proof}

\begin{definition}
\label{defn:compact} A homogeneous logical structure $(M,\cR)$ is
\emph{compact} if for each $n$, set $S$ of relations of the form
$R(x_1,\dots,x_n,a^R_1,\dots,a^R_m)$, where $R\in\cR$,
$x_1,\dots,x_n$ are fixed and $a^R_1,\dots,a^R_m\in M$ can vary
with $R$; if $|S|<|M|$ and every finite subset of $S$ is realized
in $(M,\cR)$, then $S$ is realized in $(M,\cR)$.
\end{definition}

\begin{remark}
\label{rem:defn_compact} Let $M$ be a saturated model of a
first-order theory. Let $\cR$ be the collection of definable
relations on $M$. Then $(M,\cR)$ is a compact homogeneous logical
structure. Conversely, if $(M,\cR)$ is a compact homogeneous
logical structure, consider $M$ in an expanded language with a
predicate symbol for every relation in $\cR$. In this language $M$
is a saturated model.
\end{remark}

\begin{remark}
Let $M$ be a homogeneous model of a first-order theory. Let $\cR$
be the collection of first-order definable relations on $M$. Then
$(M,\cR)$ is a homogeneous logical structure. The study of {\em
stability}  for such structures was initiated by Shelah in
\cite{Sh:4}. Our context is formally more general (as we allow
$\cR$ to stand for more general relations), but we have phrased
our definitions so that the existing stability machinery holds in
our context with obvious minor modifications. (See Section 5 on
Stability for details.)
\end{remark}

\section{Dividing and simplicity in a homogeneous logical structure}
\label{sec:div}

Throughout this section $(M,\cR)$ is an  $s\lambda-$homogeneous
logical structure such that
\begin{description}
\item[$(\bold{P})$] for some infinite cardinals
$\pi\leq\pi'\leq\lambda$ and every $\cR-$type $p(\vv)$ over
$A\subset M$, $|A|<\pi$ and $|\vv|<\pi$, if $X$ is a sequence of
realizations of $p$ of length $\geq\pi'$, then there is a sequence
$\{\,\bb_i:\:i<\omega\,\}$ $A-$indiscernible such that
$tp(\bb_0,\dots,\bb_n/A)$ is realized by an increasing sequence in
$X$, for each $n<\omega$.
\end{description}

\noindent This convention may be restated in important definitions
and results for clarity.

Indiscernible sequences play an integral role in this treatment of
dividing theory. Indeed, even the definition of dividing involves
indiscernibles. For sufficiently large $\lambda$, $(\bold{P})$
holds for any $s\lambda-$homogen-eous logical structure by
Lemma~\ref{lem:exist_indisc}.

As stated in the introduction when studying the models of a
first-order theory it is common to restrict attention to a fixed
universal domain $N$. When working in a universal domain a common
fact of model theory such as ``a consistent type over a subset of
a model can be realized in some other model'' is replaced by ``a
consistent type over a subset of $N$ of cardinality $<|N|$ is
realized in $N$. That is, only types over subsets of cardinality
$<|N|$ are studied. This restriction is realized in a convention
that the terms ``set'' and ``model'' only refer to objects of size
$<|N|$.

In analogue to the first-order conventions, given $(M,\cR)$
satisfying $(\bold{P})$, the term ``set'' will refer to a subset
of $M$ of cardinality $<\pi$. By extension the term $\cR-$type
will only apply to an $\cR-$type in $<\pi$ variables over a set of
cardinality $<\pi$. We may restate the restriction ``of
cardinality $<\pi$'' for clarity in a context where $(\bold{P})$
is being explicitly used.

\subsection{Main definitions}
\label{subsec:main_defn}

The principal concepts in this paper are ``an $\cR-$type $p$
divides over $A\subset M$'' and ``$M$ is simple''.

\begin{definition}
\label{defn:divide} Given an  $s\lambda-$homogeneous logical
structure $(M,\cR)$ satisfying $(\bold{P})$, an $\cR-$type
$p(\vv,\bb)$ over $\bb$ \emph{divides over $A\subset M$}, if there
is an infinite $A-$indiscernible sequence $\{\,\bb_i:\:i\in
X\,\}$,  with $tp(\bb_0/A)=tp(\bb/A)$, such that $\bigcup_{i\in
X}p(\vv,\bb_i)$ is inconsistent.
\end{definition}

\begin{remark}\label{rem:defn_div}
(i) Suppose $p(\vv,\bb)$ divides over $A$. Since $X$ in the
definition is infinite, Lemma~\ref{lem:large_indisc} implies that
$tp(\bb/A)$ is large. Consequently, when $tp(\dd/A)$ is small any
$\cR-$type $q(\vv,\dd)$ does not divide over $A$.

(ii) The inconsistency of $\bigcup_{i\in X}p(\vv,\bb_i)$ in the
definition depends on $tp(\{\,\bb_i:\:i\in X\,\}/A)$. The
definition does not exclude the possibility of an infinite
$A-$indiscernible sequence $\{\,\cc_i:\:i\in Y\,\}$ having the
same diagram over $A$ as $\{\,\bb_i:\:i\in X\,\}$, where
$\bigcup_{i\in Y}p(\vv,\cc_i)$ consistent.

(iii) $tp(a/b)$ divides over $C$ if and only if $tp(a/bC)$ divides
over $C$.

(iv)  If $tp(a/b)$ divides over $A$ and $C\subset A$, then
$tp(a/b)$ divides over $C$.
\end{remark}

The following basic properties of dividing do not require any
additional properties of $M$, they are properties of dividing
itself.

\begin{lemma}
\label{lem:div_prop1} Let $p(\vv,\bb)$ be an $\cR-$type over $\bb$
and $A\subset M$. The following are equivalent.
\begin{enumerate}[(1)]
\item $p(\vv,\bb)$ does not divide over $A$. \item For any
infinite $A-$indiscernible sequence $I$ there is  $J$ an infinite
indiscernible sequence over $A$ and an $\aa'$ realizing
$p(\vv,\bb')$, for $\bb'\in J$, such that $J$ has the same type
diagram over $A$ as $I$ and $J$ is indiscernible over
$A\cup\{\aa'\}$. \item For each infinite indiscernible $I$ over
$A$ with $\bb\in I$ there is an $\aa'$ realizing $p(\vv,\bb)$ such
that $I$ is indiscernible over $A\cup\aa'$.
\end{enumerate}
\end{lemma}

\begin{proof}
(2)$\;\Longrightarrow\;$(1) follows from the definition of
dividing and (3)$\;\Longrightarrow\;$(2) is trivial. To prove that
(1) implies (2) let $I=\{\,\bb_i:\:i\in X\,\}$ be an infinite
$A-$indiscernible sequence with $tp(\bb_i/A)=tp(\bb/A)$. By
Lemma~\ref{lem:large_indisc}, there is an ordered set $X'$
extending $X$ with $\pi'\leq |X'|<\lambda$ (see $(\bold{P})$ at
the beginning of the section for the definition of $\pi'$) such
that $I'=\{\,\bb_i:\:i\in X'\,\}$ is $A-$indiscernible. Since
$p(\vv,\bb)$ does not divide over $A$ there is $\aa$ in $ M$
realizing $\bigcup_{i\in X'}p(\vv,\bb_i)$.
Lemma~\ref{lem:exist_indisc} gives the existence of
$\{\,\dd_i':\:i<\omega\,\}$ indiscernible over $A\cup\aa$ such
that for $n<\omega$ there exists $i_0<\dots<i_n$ in $X'$
satisfying
\[ tp(\dd_0',\dots,\dd_n'/A)=tp(\bb_{i_0},\dots,\bb_{i_n}/A).\]
Since $p(\aa,\dd_i')$ for each $i$ we've proved (1) implies (2).

For (2) implies (3), let $I=\{\,\bb_i:\:i\in X\,\}$ and
$\{\,\dd_i':\:i<\omega\,\}$ indiscernible over $A\cup\aa$ with the
same diagram over $A$ as $I$ such that $p(\aa,\dd'_i)$, for
$i<\omega$. By Lemma~\ref{lem:large_indisc}(ii) there is a
sequence $\{\,\dd_i:\:i\in X\,\}$ indiscernible over $A\cup\aa$
indexed by $X$ with the same diagram over $A\cup\aa$ as
$\{\,\dd_i':\:i<\omega\,\}$. Now $I$ and $\{\,\dd_i:\:i\in X\,\}$
are both indiscernibles indexed by $X$ with the same diagram over
$A$. Thus, there is an automorphism $f$ of $M$ fixing $A$ and
taking $\dd_i$ to $\bb_i$, for $i\in X$. Then $f(\aa)=\aa'$
realizes $\bigcup_{i\in X}p(\vv,\bb_i)$ and $I$ is indiscernible
over $A\cup\aa'$, proving the lemma.
\end{proof}

\begin{proposition}[Pairs Lemma]
\label{prop:left_trans} If $tp(\aa/A\cup\bb)$ does not divide over
$A$ and $tp(\cc/A\cup\bb\aa)$ does not divide over $A\cup\aa$,
then $tp(\aa\cc/A\cup\bb)$ does not divide over $A$.
\end{proposition}

\begin{proof}
Let $I$ be any infinite indiscernible sequence over $A$ containing
$\bb$. By the preceding lemma there is $\aa_0$ in $M$ realizing
$tp(\aa/A\cup\bb)$ such that $I$ is indiscernible over
$A\cup\aa_0$. Let $f$ be an automorphism fixing $A\cup\bb$ and
sending $\aa$ to $\aa_0$. Then, $tp(f(\cc)/A\cup\bb\aa_0)$ does
not divide over $A\cup\aa_0$, hence there is $\cc_0$ realizing
$tp(f(\cc)/A\cup\bb\aa_0)$ such that $I$ is indiscernible over
$A\cup\aa_0\cc_0$. Since $\aa_0\cc_0$ realizes
$tp(\aa\cc/A\cup\bb)$ and $\bb\in I$ the proposition is proved.
\end{proof}

Our goal is to find minimal properties of $M$ on which dividing
defines a symmetric and transitive dependence relation, ultimately
leading to a dimension theory. Defining this property, simplicity,
will take a couple of preliminary notions.

\begin{definition}
\label{defn:free} Given an infinite cardinal $\chi$ and  $A$, $B$,
$C$ subsets of $M$, \emph{$A$ is $\chi-$free from $B$ over $C$} if
for all sequences $\aa$ from $A$ and $\bb$ from $B\cup C$ with
$|\aa|,\;|\bb|<\chi$, $tp(\aa/\bb)$ does not divide over $C$.
\end{definition}

\begin{remark}
\label{rem:defn_free} When $(M,\cR)$ is a universal domain of a
first-order theory $\aleph_0-$freeness agrees with the forking
independence relation mentioned in the paper's introduction. This
is slightly inaccurate since ``$tp(\aa/\bb)$ does not divide over
$C$'' is equivalent to ``$tp(\aa/\bb)$ does not fork over $C$''
only in a simple theory.
\end{remark}

In an arbitrary model $\kappa-$freeness depends  on types in
$<\kappa$ variables over sequences of length $<\kappa$, which may
be infinite if $\kappa$ is uncountable.

\begin{definition}
\label{defn:character} Let $\mathcal{F}$ denote the
$\chi-$freeness relation in $(M,\cR)$. The \emph{character} of
$\mathcal{F}$ is the least cardinal $\mu$ such that for all sets
$A$, $B$, $C$ in $M$, if $A$ is $\mu-$free from $B$ over $C$, then
$A$ is $\chi-$free from $B$ over $C$. $\mathcal{F}$ has
\emph{finite character} when the character is $\aleph_0$.
\end{definition}

\begin{remark}
(i) In most natural instances $\kappa-$freeness has finite
character. In models in which $\kappa-$freeness does not have
finite character we have to assume outright additional properties
of $\kappa-$freeness that lead to a notion of freeness that is
symmetric, transitive and has type amalgamation. Moreover, the
resulting notion of freeness is so esoteric that a rich theory of
dependence is unlikely. For these reasons we include finite
character in the definition of simplicity in this paper. With the
proper assumptions on $\kappa-$freeness replacing finite character
the same proofs used here work more generally. It should be noted
that finite character holds when $\kappa=\aleph_0$ and it is in
this context in which we can expect the most powerful tools of
geometrical stability theory to generalize.

(ii) When $(M,\cR)$ is compact and $\chi$ is an infinite cardinal,
$\chi-$freeness has finite character. To prove this suppose $\aa$
is a sequence of length $<\chi$, $B\subset A$ and $\aa$ is not
$\chi-$free from $A$ over $B$. Let $\cc$ be a sequence of length
$<\chi$ such that $p(\xx,\cc)=tp(\aa/\cc)$ divides over $B$. Let
$I=\{\,\cc_i:\:i\in X\,\}$ be an infinite $B-$indiscernible
sequence with $\cc_0=\cc$ such that $\bigcup_{i\in X}p(\xx,\cc_i)$
is inconsistent. By compactness, for each $i\in X$ there is a
finite $\dd_i\subset\cc_i$ such that $\bigcup_{i\in
X}p(\xx,\dd_i)$, and without loss of generality,
$J=\{\,\dd_i:\:i\in X\,\}$ is indiscernible. Thus, $tp(\aa/\dd_0)$
divides over $B$; i.e., $\aa$ is not $\aleph_0-$free from $A$ over
$B$.
\end{remark}

\begin{definition}
\label{defn:ext_base} An infinite cardinal $\kappa$ is given. A
set $A\subset M$ is a \emph{$\kappa-$extension base} if for any
sequence $\aa$ of length $<\kappa$ with $tp(\aa/A)$ large, and
$B\subset A$ such that $\aa$ is $\kappa-$free from $A$ over $B$,
and any set $C$, $A\subset C\subset M$, there is a $\cc$ realizing
$tp(\aa/A)$ such that $\cc$ is $\kappa-$free from $C$ over $B$.
\end{definition}

\begin{remark}
\label{rem:defn_ext_base} Colloquially speaking, over an extension
base  any large type has $\kappa-$free extensions over any larger
set.
\end{remark}

\begin{definition}
\label{defn:simple} Let $(M,\cR)$ be $s\lambda-$homogeneous,
$\kappa$ an infinite cardinal.

(i) $(M,\cR)$ is \emph{almost $\kappa-$simple} if it satisfies
\begin{enumerate}[(1)]
\item (Finite Character) $\kappa-$freeness has finite character.
\item  (Bounded Dividing Property) For any sequence $\aa$ of
length $<\kappa$ and set $A$, there is $B\subset A$, $|B|<\kappa$,
such that $\aa$ is $\kappa-$free from $A$ over $B$. \item (Free
Extension Property) Given a set $A\subset M$, $|A|<\pi$, there is
a $\kappa-$extension base $A'$, $A\subset A'\subset M$ and
$|A'|<\pi$.
\end{enumerate}

(ii) $(M,\cR)$ is \emph{$\kappa-$simple} if it is almost
$\kappa-$simple and every set $A\subset M$, $|A|<\pi$, is a
$\kappa-$extension base.

(iii) $(M,\cR)$ is \emph{almost simple}  (\emph{simple}) if it is
almost $\kappa-$simple ($\kappa-$simple) for some infinite
$\kappa$, with $\kappa(M)$ denoting the least such cardinal.

(iv) When $(M,\cR)$ is \emph{almost $\aleph_0-$simple}
($\aleph_0-$simple) it is also called \emph{almost supersimple}
(\emph{supersimple}).
\end{definition}

\begin{remark}
\label{rem:defn_simple} (i) First, let's compare the definition of
$\kappa-$simple with the ordinary definition of a simple theory.
To distinguish from the term defined above, \emph{classically
simple} will be used for the concept normally applied to a
first-order theory. Classically simple is defined as follows. Let
$(M,\cR)$ be a universal domain of a first-order theory viewed as
a logical structure; i.e., $\cR$ is the class of definable
relations and $(M,\cR)$ is compact. A relation $R(\xx,\aa)$
\emph{forks over $A$} if there are
$R_0(\xx,\yy_0),\dots,R_n(\xx,\yy_n)\in\cR$ (for some $n$) and
$\bb_0,\dots,\bb_n$ such that (1)
$R(\xx,\aa)\;\longrightarrow\;\bigvee_{i\leq n}R(\xx,\bb_i)$ and
(2) $R(\xx,\bb_i)$ divides over $A$ for all $i$. Then, $(M,\cR)$
is \emph{classically simple} if there is a cardinal $\kappa$ such
that for every finite sequence $\aa$ and set $A$ there is
$B\subset A$, $|B|<\kappa$, and every $R(\xx,\bb)\in tp(\aa/A)$
does not fork over $B$. Note: by the compactness of $(M,\cR)$,
$R(\xx,\aa)$ does not fork over $A$ if and only if for any set
$B\supset A$ there is a $\cc$ satisfying $R(\xx,\aa)$ such that no
$S(\xx,\bb)\in tp(\cc/B)$ divides over $B$, equivalently, $\cc$ is
$\aleph_0-$free from $B$ over $A$. By compactness $\cc$ is
$\kappa-$free from $B$ over $A$. So, the existence of free
extensions is built into the definition of not forking. It is
routine to show that a  classically simple $(M,\cR)$ is
$\kappa-$simple as defined in Definition~\ref{defn:simple}.

(ii) When $\kappa=\aleph_0$, $\kappa-$freeness has finite
character vacuously. Thus, a definition of supersimple could be
restated without this condition.

(iii) The authors' experience with dividing outside of the
first-order context suggests that there may be important examples
of models that are almost $\kappa-$simple but not $\kappa-$simple.
While little will be done with the concept here we feel it is
worthwhile examining the most fundamental properties of dividing
under this hypothesis as well as under $\kappa-$simplicity.
\end{remark}

\begin{notation}
When $(M,\cR)$ is almost $\kappa-$simple we drop $\kappa$ from the
term $\kappa-$free and simply say \emph{free} or
\emph{independent}. Remember that part of the definition of almost
$\kappa-$simple is the assumption that $\kappa-$freeness has
finite character. So, in this setting, $\aleph_0-$free implies
free.

If $p=tp(\aa/A)$ and $\aa$ is free from $B$ over $A$ then we say
$q=tp(\aa/A\cup B)$ is a \emph{free extension of $p$} or $q$ is
\emph{free from $B$ over $A$}.

If $(M,\cR)$ is almost $\kappa-$simple, $\indep{A}{C}{B}$ denotes
$A$ is free from $B$ over $C$.
\end{notation}

\begin{notation}
If $(M,\cR)$ is $\kappa-$simple letters $a,\,b,\,c,\ldots$ denote
sequences of length $<\kappa$ from $M$ and $x,\,y,\,v,\ldots$
corresponding sequences of variables.
\end{notation}

Our goal is to prove that $\kappa-$freeness is a well-behaved
dependence relation when $(M,\cR)$ is $\kappa-$simple. Symmetry
and transitivity are two critical properties. These are
straightforward consequences of
Proposition~\ref{prop:witness_divide}. To prove that proposition
we need  to prove that  the class of free extensions of a large
type is sufficiently rich. The relevant definitions and lemmas
follow.

Given an ordered set $X$ and a set $I=\{\,a_i:\:i\in X\,\}$,
$a_{<i}$ denotes $\{\,a_j:\:j<i\,\}$. If $Y,\,Z\subset X$, then
$Z<Y$ means $z<y$, for all $z\in Z$ and for all $y\in Y$. $a_Z$
denotes $\{\,a_i:\:i\in Z\,\}$.

\begin{definition}
\label{defn:morley} Let $p=tp(\aa/A)$, $\chi$ an infinite cardinal
and $B\subset A$ such that $\aa$ is $\chi-$free from $A$ over $B$.
Let $X$ be an ordered set. A sequence $\{\,\aa_i:\:i\in X\,\}$ is
a \emph{$\chi-$Morley sequence in $p$ over $B$} if
\begin{itemize}
\item each $\aa_i$, $i\in X$, realizes $p$, \item
$\{\,\aa_i:\:i\in X\,\}$ is indiscernible over $A$, and \item for
each $Y,\,Z\subset X$, $|Y|,\,|Z|<\chi$ and $Z<Y$, $\aa_Y$ is
$\chi-$free from $A\cup\aa_Z$ over $B$.
\end{itemize}
If $(M,\cR)$ is almost $\kappa-$simple, $I$ is called a
\emph{Morley sequence} if $I$ is a $\kappa-$Morley sequence.
\end{definition}

\begin{lemma}
\label{lem:special_morley} Let $(M,\cR)$ be
$s\lambda-$homogeneous, $B\subset A\subset M$, $|B|<\kappa$, and
$I=\{\,\aa_i:\:i\in X\,\}$ an infinite $A-$indiscernible sequence
such that for any $i\in X$, $\aa_i$ is $\aleph_0-$free from
$A\cup\aa_{<i}$ over $B$. Then $I$ is an $\aleph_0-$Morley
sequence over $B$.
\end{lemma}

\begin{proof}
Suppose to the contrary that $I$ is not an $\aleph_0-$Morley
sequence over $B$. Let $Y<Z$ be finite subsets of $X$ such that
for some finite subsets $\cc\subset\aa_Y$ and $\dd\subset\aa_Z$,
$tp(\dd/B\cup\cc)$ divides over $B$, and $|Z|$ is minimal with
this property. Let $i$ be the largest element of $Z$ and
$W=Z\setminus\{i\}$. Let $\dd\cap\aa_i=\ee$ and
$\dd\cap\aa_W=\ff$. Then, $tp(\ff/B\cup\cc)$ does not divide over
$B$ by the minimality assumption on $Z$. By the hypotheses of the
lemma, $tp(\ee/B\cup\cc\ff)$ does not divide over $B$. So, by
Pairs Lemma (Proposition~\ref{prop:left_trans}),
$tp(\ee\ff/B\cup\cc)$ does not divide over $B$. This contradicts
our assumption that $tp(\dd/B\cup\cc)$ divides over $B$ to prove
the lemma.
\end{proof}

\begin{lemma}
\label{lem:ext_imp_Morley} Let $(M,\cR)$ be almost
$\kappa-$simple, $A$ an extension base, $\aa$ a sequence of length
$<\kappa$ such that $p=tp(\aa/A)$ is large, and $B\subset A$ such
that $\aa$ is $\kappa-$free from $A$ over $B$. Let $X$ be any
infinite linear order with $|X|\leq\lambda$. Then, $M$ contains an
$A-$indiscernible sequence $I=\{\,\aa_i:\:i\in X\,\}$ which is a
$\kappa-$Morley sequence in $p$ over $B$.
\end{lemma}

\begin{proof}
First find an infinite $A-$indiscernible sequence that is
$\aleph_0-$free using Lemma~\ref{lem:exist_indisc} as follows. Let
$\alpha$ be an ordinal $<\lambda$ and suppose sequences $\aa_i$,
$i<\alpha$, have been chosen so that $\aa_i$ realizes $p$ and
$\aa_i$ is $\kappa-$free from $A\cup\aa_{<i}$ over $B$. By the
Extension Property in the definition of $\kappa-$simple there is
$\aa_\alpha$ realizing $p$ which is $\kappa-$free from
$A\cup\aa_{<\alpha}$ over $B$. Since $(M,\cR)$ satisfies
$(\bold{P})$ there is $J=\{\,\bb_i:\:i<\omega\,\}$ such that $J$
is $A-$indiscernible and for any $n<\omega$ there are
$i_1,\dots,i_n<\lambda$ with $tp(\bb_1,\dots,\bb_n/A) =
tp(\aa_{i_1},\dots,\aa_{i_n}/A)$. This latter property controls
the type diagram of $J$ and by Lemma~\ref{lem:special_morley}
guarantees that $J$ is an $\aleph_0-$Morley sequence.

Given an infinite linear order $X$, $|X|\leq\lambda$, by
Lemma~\ref{lem:large_indisc}(ii) there is a sequence
$K=\{\cc_i:\:i\in X\,\}$  indiscernible over $A$ with the same
type diagram as $J$. Thus, $K$ is also an $\aleph_0-$Morley
sequence in $p$ over $B$. Since $(M,\cR)$ is $\kappa-$simple,
$\kappa-$freeness has finite character. Thus, $K$ is also a
$\kappa-$Morley sequence, proving the lemma.
\end{proof}

\begin{remark}
The proof of the preceding lemma illustrates the problems
circumvented by assuming finite character. Being a $\kappa-$Morley
sequence depends on all subsequences of length $<\kappa$. Without
finite character the partition calculus required to get a
$\kappa-$Morley sequence could easily be independent of set theory
(depending on $\kappa$ and $|M|$). For the same reason it is
unclear that a $\kappa-$Morley sequence can be extended to a
larger $\kappa-$Morley sequence. While the definition of
simplicity could be rewritten to force the existence of enough
Morley sequences this would simply hide the complexity in a
definition. We believe the class of models satisfying the
necessary conditions (without assuming finite character) is very
thin, hence our assumption of finite character as part of
Definition~\ref{defn:simple}.
\end{remark}

The following proposition is the key to obtaining symmetry and
transitivity of freeness. It says that when $p(\xx,\bb)$ divides
over $A$, this fact can be witnessed with a Morley sequence in
$tp(\bb/A)$.

\begin{proposition}
\label{prop:witness_divide} Let $(M,\cR)$ be almost
$\kappa-$simple, $A$ an extension base, $|b|<\kappa$, and $p(x,b)$
a type over $b$ with $|x|<\kappa$. An indiscernible sequence $I$
indexed by $X$ satisfies $(\ast)$ if for any suborder $Y$ of $X$
with $|Y|<\kappa$ there is $i\in X$, $i<Y$. Suppose $b$ is free
from $A$ over $B\subset A$, $|B|<\kappa$. The following are
equivalent.
\begin{enumerate}[(1)]
\item $p(x,b)$ divides over $A$. \item There is a Morley sequence
$I$ in $tp(b/A)$ over $B$ satisfying  $(\ast)$ such that
$\bigcup_{d\in I}p(x,d)$ is inconsistent. \item For any  Morley
sequence $I$ in $tp(b/A)$ over $B$ satisfying  $(\ast)$,
$\bigcup_{d\in I}p(x,d)$ is inconsistent.
\end{enumerate}
\end{proposition}

\begin{proof}
(2)$\;\Longrightarrow\;$(1) is simply by the definition of
dividing. (3)$\;\Longrightarrow\;$(2) because such a Morley
sequence exists by Lemma~\ref{lem:ext_imp_Morley}.

We now prove (1)$\;\Longrightarrow\;$(3). Let $I=\{\,b_i:\:i\in
X\,\}$ be the given Morley sequence and suppose $p(x,b)$ divides
over $A$. Notice that $I$ is also a Morley sequence in $tp(b/B)$
over $B$ and satisfies $(\ast)$, and $p(x,b)$ divides over $B$.
So, replacing $A$ by $B$ if necessary, we may as well assume
$|A|<\kappa$.
\begin{claim}
For each $i\in X$ and $Y\subset X$ with $|Y|<\kappa$ and $i<Y$,
$p(x,b_i)$ divides over $A\cup b_Y$.
\end{claim}
Simply because $b_i$ realizes $tp(b/A)$, $p(x,b_i)$ divides over
$A$. This fact is witnessed by an infinite $A-$indiscernible
sequence $J$ containing $b_i$. Let $q(z,y)=tp(b_i,b_Y/A)$. Since
$b_Y$ is $\kappa-$free from $\{\,b_j:\:j<Y\,\}$ over $A$ and
$|A|<\kappa$, $q(b_i,y)$ does not divide over $A$. Thus, there is
a sequence $c$ such that $\bigcup_{d\in J}q(d,c)$. By
Lemma~\ref{lem:div_prop1}, we can assume that $J$ is indiscernible
over $A\cup c$. Since $tp(b_ic/A)=tp(b_ib_Y/A)$ there is an
automorphism $f$ fixing $A\cup b_i$ taking $c$ to $b_Y$. Let
$J'=f(J)$. Then, $J'$ is infinite, indiscernible over $A\cup b_Y$
and $\bigcup_{d\in J'}p(x,d)$ is inconsistent. This witnesses that
$p(x,b_i)$ divides over $A\cup b_Y$.

To continue with the proof suppose to the contrary that
$\bigcup_{d\in I}p(x,d)$ is realized by some $a$. By the Bounded
Dividing Property there is $Y\subset X$ of cardinality $<\kappa$
such that $a$ is $\kappa-$free from $A\cup I$ over $A\cup b_Y$.
Since $X$ satisfies $(\ast)$ there is $i\in X$, $i<Y$. By the
claim, $p(x,b_i)$ divides over $A\cup b_Y$. Since $a$ satisfies
$p(x,b_i)$ we contradict that $a$ is $\kappa-$free from $A\cup I$
over $A\cup b_Y$. This proves the proposition.
\end{proof}

\subsection{Small types, large types and dividing}

In a saturated model  algebraic types and algebraic closure play a
very special role. On the one hand, dividing trivializes for
algebraic types: If $tp(a/A)$ is algebraic then $tp(a/A\cup\{b\})$
and $tp(b/A\cup\{a\})$ do not divide over $A$ or all $b$. On the
other hand, the close relationship between algebraic closure and
dividing in a supersimple model (through minimal types and
canonical bases) is the basis for geometrical stability theory.

In a saturated model the algebraic types are exactly the small
types (see Definition~\ref{defn:large}). We will show here that
small types and small closure act much like algebraic types in a
simple model, and to  a lesser degree in an almost simple model.

\begin{notation}
Let $s\cl(-)$ denote the \emph{small closure operator} defined by
$\aa\in s\cl(A)$ if $tp(\aa/A)$ is small.
\end{notation}

\begin{remark}
\label{rem:small_closure} For any $s\lambda-$homogeneous
$(M,\cR)$, $s\cl(-)$ is a closure operator. That is, $s\cl(-)$
satisfies:
\begin{enumerate}[(i)]
\item $A\subset s\cl(A)$; \item $A\subset
B\;\Longrightarrow\;s\cl(A)\subset s\cl(B)$; \item
$s\cl(s\cl(A))=s\cl(A)$.
\end{enumerate}
[The proof is left to the reader.]
\end{remark}

In any saturated model, (1) if $a\in acl(A)$ then for all $b$,
neither $tp(b/A\cup\{a\})$ nor $tp(a/A\cup\{b\})$  divides over
$A$, and (2) $a\in acl(A\cup\{b\})\setminus acl(A)$ implies that
$tp(a/A\cup\{b\})$ and $tp(b/A\cup\{a\})$ both divide over $A$.
Property (1) extends to small closure directly, while
Proposition~\ref{prop:main_get_small} serves as the generalization
of (2).

\begin{lemma}
\label{lem:small_clos_div} Let $(M,\cR)$ be $s\lambda-$homogeneous
and $a\in s\cl(A)$. Then for any $b$, $tp(a/A\cup\{b\})$ and
$tp(b/A\cup\{a\})$ do not divide over $A$.
\end{lemma}

\begin{proof}
If $tp(b/A)$ is small then $tp(a/A\cup\{b\})$ does not divide over
$A$ (Remark~\ref{rem:defn_div}). Suppose $q(x,b)=tp(a/A\cup\{b\})$
is large and let $I$ be an arbitrary $A-$indiscernible sequence in
$tp(b/A)$. For any $d\in I$, any realization of $q(x,d)$ also
realizes the small type $tp(a/A)$. Using
Lemma~\ref{lem:exist_indisc} and cardinality properties we find a
$c$ realizing $tp(a/A)$ and an infinite sequence $J$ indiscernible
over $A\cup\{c\}$ with the same type diagram over $A$ as $I$, such
that $c$ realizes $\bigcup_{d\in J}q(x,d)$. By
Lemma~\ref{lem:div_prop1}, $q(x,b)$ does not divide over $A$.

Using $a\in s\cl(A)$ and Lemma~\ref{lem:large_indisc} there is no
infinite $A-$indiscernible sequence in $tp(a/A)$. Thus,
$tp(b/A\cup\{a\})$ does not divide over $A$.
\end{proof}

\begin{proposition}
\label{prop:main_get_small} Let $(M,\cR)$ be $\kappa-$simple,
$a,\,b$ and $c$ sequences of length $<\kappa$ such that $tp(a/c)$
is large and $tp(a/b,c)$ is small. Then, $tp(a/b,c)$ divides over
$c$ and $tp(b/a,c)$ divides over $c$.
\end{proposition}

The complete proof will take several lemmas. Proving that the
second of the types divides is straightforward and holds in an
almost simple model.

\begin{lemma}
\label{lem:get_small} Let $(M,\cR)$ be almost $\kappa-$simple,
 $a$, $b$ and $c$ sequences of cardinality $<\kappa$
with $tp(a/c)$ large and $tp(a/b,c)$ small. Then, $tp(b/a,c)$
divides over $c$.
\end{lemma}

\begin{proof}
This is simply a  counting argument. Let $q(x,y)=tp(a,b/c)$. Since
$tp(a/c)$ is large there is an infinite $c-$indiscernible sequence
$I$ in $tp(a/c)$. Since $q(x,d)$ is small for any $d$ there cannot
be a $d$ realizing $tp(b/A)$ such that  $I$ is indiscernible over
$c,d$ and $q(a',d)$, for $a'\in I$. By Lemma~\ref{lem:div_prop1},
$q(a,y)$ divides over $c$.
\end{proof}

Completing the proof requires a couple of lemmas. While none of
the results are especially significant, they are separated into
independent lemmas for ease of reference. The following lemma is
immediate by the definition of dividing.

\begin{lemma}
\label{lem:alg_div} If $tp(a/A)$ is large, then $x=a$ divides over
$A$.
\end{lemma}

\begin{lemma}
\label{lem:get_small_div} Let $(M,\cR)$ be almost $\kappa-$simple,
$a,\,b$ sequences of length $<\kappa$, $A$ an extension base such
that $tp(a/A)$ is large and $tp(a/B\cup\{b\})$ is small for some
$B\subset A$, $|B|<\kappa$. Then, $tp(a/B\cup\{b\})$ divides over
$A$.
\end{lemma}

\begin{proof}
Let $q(x,b)=tp(a/B\cup\{b\})$. Let $C\subset A$, $|C|<\kappa$, be
such that $b$ is free from $A$ over $C$ and $B\subset C$. Notice
that $tp(b/A)$ is large since $tp(a/A)$ is large and
$tp(a/A\cup\{b\})$ is small. Let $\mu<\lambda$ be a cardinal
$>\kappa$ and $I=\{\,b_i:\:i<\mu^*\,\}$ be a Morley sequence in
$tp(b/A)$ over $C$, where $\mu^*$ is the reverse order on $\mu$.
For $i<\mu^*$ let $\cc_i$ enumerate the set of realizations of
$q(x,b_i)$. Since $q(x,b_i)$ is small we can choose $I$ so that
$\{\,b_j:\:j>i\,\}$ is indiscernible over
$A\cup\{\,\cc_kb_k:\:k\leq i\,\}$ and free from
$A\cup\{\,\cc_kb_k:\:k\leq i\,\}$ over $C$. Given $c'\subset
\cc_i$ an arbitrary sequence of length $<\kappa$, $c'$ is free
from $A\cup\{\,\cc_kb_k:\:k<i\,\}$ over $B\cup\{b_i\}$ since
$tp(c'/B\cup\{b_i\})$ is small (Lemma~\ref{lem:small_clos_div}).
By Pairs Lemma (Proposition~\ref{prop:left_trans}), $c'b_i$ is
free from $A\cup\{\,\cc_kb_k:\:k<i\,\}$ over $C$, hence,
$\cc_ib_i$ is free from $A\cup\{\,\cc_kb_k:\:k<i\,\}$ over $C$. In
particular, $\cc_i$  is free from $\{\,\cc_k:\:k<i\,\}$ over $C$.
By Lemma~\ref{lem:alg_div} the sequence $\cc_i$ is disjoint from
$\bigcup_{k<i}\cc_k$. We conclude that $q(x,b_i)\cup q(x,b_j)$ is
inconsistent for $i<j<\lambda^*$. Since $I$ is a Morley sequence,
$q(x,b)$ divides over $A$ (Proposition~\ref{prop:witness_divide}).
\end{proof}

\begin{remark}
\label{rem:on_smallness} Let $(M,\cR)$ be almost $\kappa-$simple.
Properties of the small closure operator significantly affect
dividing and the overall structure of the model. The preceding
lemma and Remark~\ref{rem:defn_div}(i) show there are
relationships between freeness and small closure, but the picture
may be very complicated. To organize this discussion, let $\chi_s$
denote the \emph{character of smallness}; i.e., the least cardinal
$\chi$ such that if $|\aa|<\kappa$ and $tp(\aa/A)$ is small, then
there is $B\subset A$, $|B|<\chi$, such that $tp(\aa/B)$ is small.
If $\chi_s\leq\kappa$, then $tp(\aa/A)$ large and $\aa$ free from
$B$ over $A$ implies that $tp(\aa/B)$ is large. When
$\chi_s>\kappa$, which can't be ruled out in general, certain
results require restricting explicitly to large types or small
types. Better results are possible in stable theories.
\end{remark}

\subsection{Symmetry and Transitivity}
\label{subsec:symm}

We begin with versions of symmetry and transitivity for sequences
of length $<\kappa$ and generalize subsequently.

\begin{lemma}[Local Symmetry Lemma]
\label{lem:local_symm} Let $(M,\cR)$ be  $\kappa-$simple, $a,\,b$
and $c$ sequences of length $<\kappa$ such that $tp(a/bc)$ does
not divide over $c$. Then $tp(b/ac)$ does not divide over $c$.
\end{lemma}

\begin{proof}
If $tp(a/c)$ is small then $tp(b/ac)$ does not divide over $c$
(Lemma~\ref{lem:small_clos_div}). Suppose $tp(a/c)$ is large.
Since $tp(a/bc)$ does not divide over $c$, $tp(a/bc)$ is large
(Proposition~\ref{prop:main_get_small}). Thus, there is a Morley
sequence $I=\{\,a_i:\:i\in\kappa^*\,\}$ in
 $tp(a/bc)$ over $c$, where $\kappa^*$ is the reverse order on $\kappa$.
Let $q(x,a)=tp(b/ac)$. Since $I$ is a sequence of indiscernibles
over $bc$ in $tp(a/bc)$, $b$ realizes
$\bigcup_{i\in\kappa^*}q(x,a_i)$. By
Proposition~\ref{prop:witness_divide}, $q(x,a)$ does not divide
over $c$, proving the lemma.
\end{proof}

Following is the version of transitivity that holds for sequences
of length $<\kappa$ in a  simple theory. As this lemma
illustrates, transitivity is simply a combination of symmetry and
Pairs Lemma.

\begin{lemma}[Local Transitivity]
\label{lem:local_trans} Let $(M,\cR)$ be  $\kappa-$simple,
$c\subset b\subset a$ sequences of length $<\kappa$ and $d$,
$|d|<\kappa$, such that $tp(d/a)$ does not divide over $b$ and
$tp(d/b)$ does not divide over $c$. Then $tp(d/a)$ does not divide
over $c$.
\end{lemma}

\begin{proof}
By Local Symmetry $tp(b/cd)$ does not divide over $c$ and
$tp(a/bd)$ does not divide over $b$. Applying Pairs Lemma and
$c\subset b\subset a$, $tp(ab/cd)$ does not divide over $c$;
hence, $tp(a/cd)$ does not divide over $c$. Thus, $tp(d/a)$ does
not divide over $c$, again by Local Symmetry.
\end{proof}

Extending these results about dividing on sequences of length
$<\kappa$ to properties of freeness on arbitrary sets requires the
following.

\begin{lemma}[Weak Transitivity]
\label{lem:very_wk_trans} Let $(M,\cR)$ be  $\kappa-$simple,
$a,\,b$ sequences of length $<\kappa$ and $A$ a set such that $a$
is free from $b$ over $A$. Suppose $ab$ is free from $A$ over
$c\subset A$, $|c|<\kappa$. Then $a$ is free from $A\cup\{b\}$
over $c$.
\end{lemma}

\begin{proof}
To begin we claim that $tp(a/cb)$ does not divide over $c$.
Suppose, to the contrary, that $tp(a/cb)$ divides over $c$. Let
$q(x,y)=tp(ab/c)$. Since $ab$ is free from $A$ over $c$, $b$ is
free from $A$ over $c$. Since $tp(b/c)$ must be  large (or $a$
would be free from $b$ over $c$)
 there is an infinite $I$ that is a Morley sequence over $c$ in
 $tp(b/A)$ and indexed by $\kappa^*$.
Using that $a$ is free from $b$ over $A$ and $I$ is indiscernible
over $A$, $\bigcup_{d\in I}q(x,d)$ is consistent. Since $I$ is a
Morley sequence over $c$, this contradicts that $q(x,b)$ divides
over $c$ and Proposition~\ref{prop:witness_divide}.

Let $d\subset A$ be a sequence of length $<\kappa$. By symmetry
and the choice of $c$, $tp(d/cab)$ does not divide over $cb$. (In
fact, the type doesn't divide over $c$, but that is more than we
need.) Combining this fact and the claim with Pairs Lemma,
$tp(db/ca)$ does not divide over $c$. By symmetry, $a$ is free
from $bd$ over $c$. Thus, $a$ is free from $A\cup\{b\}$ over $c$.
\end{proof}

\begin{theorem}[Symmetry Lemma]
\label{thm:symm_lemma} Let $(M,\cR)$ be  $\kappa-$simple, $A$, $B$
and $C$ such that $A$ is free from $B$ over $C$. Then $B$ is free
from $A$ over $C$.
\end{theorem}

\begin{proof}
Let $a\subset A$, $b\subset B$ and $c_0\subset C$ be arbitrary
sequences of length $<\kappa$. Let $c_0\subset c\subset C$,
$|c|<\kappa$, be such that $ab$ is free from $C$ over $c$. Since
$a$ is free from $b$ over $C$, Lemma~\ref{lem:very_wk_trans}
implies that $tp(a/bc)$ does not divide over $c$. By Local
Symmetry (Lemma~\ref{lem:local_symm}) $tp(b/ac)$ does not divide
over $c$. Thus $tp(b/a)$ does not divide over $c$ and $tp(b/a)$
does not divide over $C$ (see Remark~\ref{rem:defn_div}) Thus, $b$
is free from $a$ over $C$. This proves the theorem.
\end{proof}

\begin{corollary}[Transitivity]
\label{cor:trans} Let $(M,\cR)$ be  $\kappa-$simple, $C\subset
B\subset A$ and $D$ such that $D$ is free from $A$ over $B$ and
$D$ is free from $B$ over $C$. Then, $D$ is free from $A$ over
$C$.
\end{corollary}

\begin{proof}
Let $d\subset D$ and $a\subset A$ be arbitrary sequences of length
$<\kappa$. Let $b\subset B$, $|b|<\kappa$, be  so that $da$ is
free from $B$ over $b$, and let $c\subset C$, $|c|<\kappa$, be
such that $db$ is free from $C$ over $c$. By Weak Transitivity
(Lemma~\ref{lem:very_wk_trans}) $d$ is free from $a$ over $b$ and
$d$ is free from $b$ over $c$. By Local Transitivity, $d$ is free
from $a$ over $c$. Thus, $d$ is free from $a$ over $C$, proving
the corollary.
\end{proof}

The following applications of symmetry and transitivity will be
used later.

\begin{corollary}
\label{cor:morley2} Let $(M,\cR)$ be  $\kappa-$simple and $I$ an
infinite Morley sequence over $A$ where $a\in I$ has length
$<\kappa$.

(i) Given $a\in I$, let $A_0\subset A$ be a set of cardinality
$<\kappa$ such that $a$ is free from $A$ over $A_0$. Then, $I$ is
a Morley sequence in $tp(a/A)$ over $A_0$.

(ii)  If $I$ is indiscernible over $B\supset A$ then $I$ is a
Morley sequence in $tp(a/B)$ over $A$.
\end{corollary}

\begin{proof}
Let $I=\{\,a_i:\:i\in X\,\}$. (i) By Transitivity, for any $i\in
X$, $a_i$ is free from $A\cup a_{<i}$over $A_0$. By
Lemma~\ref{lem:special_morley}, $I$ is a Morley sequence in
$tp(a/A)$ over $A_0$.

(ii) Without loss of generality,  $X$ is isomorphic to $\kappa^*$.
Let $b\subset B$ and  $c\subset A$  be sequences of length
$<\kappa$, $Y\subset X$, $|Y|<\kappa$, $d=I_Y$ and
$q(x,y,z)=tp(b,d,a_i/c)$, for $i\in X$, $Y<i$. Let
$J=\{\,a_i:\:i\in X,\;Y<i\,\}$. Since $q(b,d,a)$ holds for all
$a\in J$, $q(xy,a)$ does not divide over $c$ by
Proposition~\ref{prop:witness_divide}. By symmetry, $tp(a_i/bdc)$
does not divide over $c$, for $Y<i$. Thus, $I$ is a Morley
sequence in $tp(a/B)$ over $A$ (using
Lemma~\ref{lem:special_morley}).
\end{proof}

In an almost simple model symmetry holds for some special sets.

\begin{lemma}[Almost Symmetry Lemma]
\label{lem:symm} Let $(M,\cR)$ be almost $\kappa-$simple, \newline
$A\subset B$ with
 $B$ an extension base and $b$ such that $tp(b/B)$ is
large and $b$ is free from $B$ over $A$. Then, $B$ is free from
$b$ over $A$.
\end{lemma}

\begin{proof}
The hypotheses of the lemma and Lemma~\ref{lem:ext_imp_Morley}
yield a Morley sequence $I=\{\,a_i:\:i\in\kappa^*\,\}$ in
$tp(b/B)$ over $A$, where $\kappa^*$ is the reverse order on
$\kappa$. A fortiori, $I$ is a Morley sequence over $A$. Let
$q(x,b)$ be a type over $c\cup b$ for some $c\subset A$,
$|c|<\kappa$, satisfied by some $d\in B$, $|d|<\kappa$. Since $I$
is indiscernible in $tp(b/B)$, $q(d,a_i)$ holds for all
$i\in\kappa^*$. Since $I$ is a Morley sequence over $A$, $q(x,b)$
does not divide over $A$ by Proposition~\ref{prop:witness_divide}.
Thus, $B$ is free from $b$ over $A$.
\end{proof}

More general forms of transitivity in arbitrary almost simple
theories come down to applying Pairs Lemma and the form of the
Symmetry Lemma that holds in that context.

\begin{definition}
\label{defn:indep} Let $(M,\cR)$ be simple. A set of elements $I$
is \emph{independent} over $A$ if for all $a\in I$, $a$ is free
from $I\setminus\{a\}$ over $A$.
\end{definition}

\begin{remark}\label{rem:indep}
By Symmetry and Transitivity, if $I$ is a Morley sequence over $A$
then  $I$ is independent over $A$.
\end{remark}

\section{Type amalgamation}
\label{sec:tp_amalg}

The main result in the section is colloquially known as ``type
amalgamation'' and stated as Theorem~\ref{thm:tp_amal}. This
result lends insight into the question of when distinct free
extensions of a type have a common free extension. This is
critical to applying the freeness relation to  induce geometrical
structure properties. In particular, properties of the parallelism
relation (see Section~\ref{subsec:paral}) depend heavily on type
amalgamation. The theorem involves the concept of a Lascar strong
type, developed in the next subsection.

In the first-order setting the Type Amalgamation Theorem was
originally called the Independence Theorem. A rudimentary version
was found in \cite{sh:simp}. The theorem here generalizes the
first-order version proved by Kim and Pillay \cite{ki.pi:simple}.

\subsection{Lascar strong types}
\label{subsec:lascar}

In this subsection $(M,\cR)$ is an arbitrary
$s\lambda-$homogeneous model satisfying $(\bold{P})$. As usual,
all sets referenced are considered to be subsets of $M$ of
cardinality $<\pi$. Dividing theory is not used here. The notion
of Lascar strong type was introduced by Lascar in \cite{la:lstp}.

We let $SE^\mu(A)$ be the set of $A$-invariant equivalence
relations on $M^\mu$ with a bounded ($<\lambda$) number of
equivalence classes. Let $SE(A) = \bigcup_\mu SE^\mu(A)$.

\begin{definition}\label{defn:lstp} Tuples $a$, $b$ of the same length
have the same {\em Lascar strong type} over $A\subset M$, written
$lstp(a/A)= lstp(b/A)$, if $E(a,b)$ whenever $E \in
SE^{\ell(a)}(A)$.
\end{definition}

\begin{lemma}
\label{lem:small_lstp} If $tp(a/A)$ is small and
$lstp(b/A)=lstp(a/A)$ then $b=a$.
\end{lemma}

\begin{proof}
Let $p=tp(a/A)$. The equivalence relation defined by
\[(p(x)\longleftrightarrow p(y))\wedge(p(x)\longrightarrow x=y)\]
is $A-$invariant with a bounded number of classes (since $p$ is
small).
\end{proof}

\begin{lemma} \label{l:first}
If $I$ is an infinite indiscernible sequence over $A$ then
$E(a,b)$, for any $E \in SE^{\ell(a)}(A)$ and $a,b \in I$.
\end{lemma}

\begin{proof}
By Lemma~\ref{lem:large_indisc} and the $A-$invariance of $E$, $I$
can be extended to any length $\leq\lambda$. If $\neg E(a,b)$ then
$\neg E(c,d)$ for any $c \not = d$ in $I$. Hence, there are
unboundedly many equivalence classes.
\end{proof}

The (proof of the) next lemma shows that equality of Lascar strong
types over $A$ is the finest equivalence relation in $SE(A)$.
Thus, there are fewer than $\lambda$ Lascar strong types over $A$.

\begin{lemma}\label{l:main}
For tuples $a$, $b$ of the same length and a set $A$ with
$tp(a/A)$ and $tp(b/A)$ large, the following are equivalent:
\begin{enumerate}[(1)]
\item There exists $n < \omega$ and $a=a_0$, $a_1$, $a_2, \dots,
a_n=b$ such that for each $i < n$ there exists an infinite
$A$-indiscernible sequence containing $a_i$ and $a_{i+1}$; \item
$lstp(a/A)=lstp(b/A)$.
\end{enumerate}
\end{lemma}

\begin{proof} (1) implies (2) follows from the previous lemma and
transitivity of equivalence. For (2) implies (1), call $E$ the
equivalence relation defined by (1). Notice that $E$ is
$A$-invariant. Suppose $E$ had unboundedly many equivalence
classes and let $\{\, a_i :\: i < \mu \,\}$ be inequivalent
elements for some suitably large $\mu\geq\lambda$. By
Lemma~\ref{lem:exist_indisc}, there exists $\{\, d_i:\: i <
\omega\,\} $ indiscernible over $A$, such that $tp(d_0, d_1/A) =
tp(a_{j_0}, a_{j_1}/A)$ for some $j_0 < j_1 < \mu$. By the
definition of $E$, we have $E(d_0,d_1)$. Hence, $E(a_{j_0},
a_{j_1})$, by $A$-invariance, contradicting the choice of $\{\,
a_i :\: i < \mu \,\}$. Thus, $E \in SE^{\ell(a)}(A)$. Therefore,
if $lstp(a/A) = lstp(b/A)$, then $E(a,b)$ holds so that (1) holds.
\end{proof}

\begin{definition} Let $\Saut{A}$ be the set
of $f \in \aut{A}$ such that  for each $a \in M$, $lstp(f(a)/A) =
lstp(a/A)$.
\end{definition}

$\Saut A$ is a group, called  the \emph{group of strong
automorphisms over $A$}. Furthermore, $\Saut A$ is a normal
subgroup of $\aut{A}$: if $f \in \Saut A$ and $g \in \aut{A}$,
then $E(g(a), fg(a))$ by definition of $f$, hence $E(a,
g^{-1}fg(a))$, by $A$-invariance.

\begin{lemma} \label{l:strong} The following conditions are equivalent:
\begin{enumerate}[(1)]
\item $lstp(a/A)= lstp(b/A)$ \item There exists $ f\in \Saut A$
such that $f(a)=b$.
\end{enumerate}
\end{lemma}

\begin{proof}
(2) implies (1) follows immediately from the definition. We show
(1) implies (2). Define $E(a,b)$ if there exists $f \in \Saut A$
with $f(a)=b$. This is clearly an equivalence relation since
$\Saut A$ is a group. Notice that it is also $A$-invariant since
$\Saut A$ is normal in $\aut{A}$. Hence, it is enough to show that
$E \in SE(A)$. Suppose not and let $\{\, a_i :\: i < \lambda^+\,
\}$ be a large set of $E$-inequivalent elements. Let $B$ be a
bounded set extending $A$ containing a representative of every
Lascar strong type over $A$. By the pigeonhole principle, there
exists $i < j$ such that $tp(a_i/B)= tp(a_j/B)$. Let $g \in
\aut{B}$ such that $g(a_i)=a_j$. But $g \in \Saut A$, a
contradiction.
\end{proof}

\begin{corollary}
\label{cor:lstp_over_small} If $lstp(a/A)=lstp(b/A)$ and $c\in
s\cl(A)$, then $lstp(a/A\cup\{c\})=lstp(b/A\cup\{c\})$.
\end{corollary}

\begin{proof}
This is immediate by Lemma~\ref{l:strong} and
Lemma~\ref{lem:small_lstp}.
\end{proof}

\begin{lemma}\label{lem:fatten_lstp} If $lstp(a/A)=lstp(a'/A)$ and $b$ is
given, then there is $b' \in M$ such that
$lstp(ab/A)=lstp(a'b'/A)$.
\end{lemma}

\begin{proof} By Lemma~\ref{l:strong}, there exists $f \in \Saut{A}$ such
that $f(a)=a'$. Let $b' = f(b)$. Another application of the that
lemma shows the conclusion.
\end{proof}

In a simple model the Extension Property extends to Lascar strong
types.

\begin{lemma}[Strong Extension]
\label{lem:str_ext} Let $(M,\cR)$ be  $\kappa-$simple, $B\supset
A$ and $a$ such that $tp(a/A)$ is large. Then, there is $b$ such
that $lstp(b/A)=lstp(a/A)$ and $b$ is free from $B$ over $A$.
\end{lemma}

\begin{proof}
Since  $tp(a/A)$ is large there is an infinite Morley sequence $I$
over $A$ with $a\in I$. Using Lemma~\ref{lem:large_indisc} there
is an infinite sequence of indiscernibles $J$ such that
$I\hat{\;}J$ is indiscernible over $A$ and $J$ is indiscernible
over $B$. Then $J$ is a Morley sequence over $A$ and $b\in
J\;\Longrightarrow\; lstp(b/A)=lstp(a/A)$. By
Corollary~\ref{cor:morley2}, given $b\in J$, $J$ is a Morley
sequence in $tp(b/B)$ over $A$. In particular, $b$ is free from
$B$ over $A$.
\end{proof}

\subsection{Type amalgamation theorem}
\label{subsec:tp_amal_thm}

This subsection is devoted to the proof of

\begin{theorem}[Type Amalgamation]
\label{thm:tp_amal} Let $(M,\cR)$ be  $\kappa-$simple, $c$, and
$b_i,\,a_i$, for $i=1,2$, sequences of length $<\kappa$ such that
\begin{enumerate}[(1)]
\item $tp(b_1/cb_2)$ does not divide over $c$, \item
$lstp(a_1/c)=lstp(a_2/c)$, and \item $tp(a_i/cb_i)$ does not
divide over $c$.
\end{enumerate}
Then, there is $a$  realizing $lstp(a_i/cb_i)$, for $i=1,2$, such
that $tp(a/cb_1b_2)$ does not divide over $c$.
\end{theorem}

Throughout the subsection $(M,\cR)$ is   $\kappa-$simple. The bulk
of the proof will be found in preliminary lemmas that are actually
special cases of the theorem.
Proposition~\ref{prop:witness_divide} is the main preliminary
result here.

\begin{lemma}
\label{lem:morley_conj1} Let $p(x,b)$ be a type over $A\cup\{b\}$,
where $|b|<\kappa$, which does not divide over $A$, and $I$ an
infinite Morley sequence in $tp(b/A)$. Then, for any $b_0,\,b_1\in
I$, $p(x,b_0)\cup p(x,b_1)$ does not divide over $A$.
\end{lemma}

\begin{proof}
Without loss of generality, $I=\{\,a_i:\:i\in\kappa^*\,\}$ where
$\kappa^*$ is the reverse order on $\kappa$. Let $X$ be a suborder
of $\kappa^*$ such that $X$ is coinitial in $\kappa^*\setminus X$,
$\kappa^*\setminus X$ is coinitial in $X$ and $X$ is isomorphic to
$\kappa^*$. Let $f$ be an injective function from $X$ into
$\kappa^*\setminus X$ such that $i<f(i)$, for $i\in X$. Let
$J=\{\,a_ia_{f(i)}:\:i\in X\,\}$. It is routine to show that $J$
is a Morley sequence in $tp(d/A)$, for $d\in J$. For
$d=a_ia_{f(i)}\in J$, let $q(x,d)=p(x,a_i)\cup p(x,a_{f(i)})$.

Since $p(x,b)$ does not divide over $A$ there is a $c$ realizing
$\bigcup_{e\in I}p(x,e)$. Thus, $c$ realizes $\bigcup_{d\in
J}q(x,d)$. By Proposition~\ref{prop:witness_divide}, $q(x,d)$ does
not divide over $A$, proving the lemma.
\end{proof}

\begin{lemma}
\label{lem:all_r_morley} Let $A$ be a set of cardinality
$<\kappa$, $I=\{\,a_i:\:i\in X\,\}$ an infinite $A-$indiscernible
sequence with $|a|<\kappa$, for $a\in I$. Then there is $Y\subset
X$, $|Y|<\kappa$, such that $K=\{a_i:\:Y<i,\,i\in X\,\}$ is a
Morley sequence over $A\cup \{\,a_i:\:i\in Y\,\}$.
\end{lemma}

\begin{proof}
Let $X'=X+\{x\}$ be the order obtained by adding a single element
$x$ to the end of $X$. Let $a_x$ be a sequence in $M$ such that
$I'=\{\,a_i:\:i\in X'\,\}$ is $A-$indiscernible. By the
$\kappa-$simplicity of $M$ there is $Y\subset X$, $|Y|<\kappa$,
such that $a_x$ is $\kappa-$free from $A\cup I$ over $A\cup J$,
$J=\{\,a_i:\:i\in Y\,\}$. Let $Z=\{\,j\in X:\:Y<j\,\}$ and
$K=\{\,a_i:\:i\in Z\,\}$. If $i\in Z$, then $a_i$ and $a_x$ have
the same type over $A\cup J\cup (K\cap a_{<i})$. Thus, $a_i$ is
$\kappa-$free from $K\cap a_{<i}$ over $A\cup J$. By
Lemma~\ref{lem:special_morley} and $\kappa-$simplicity, $K$ is a
Morley sequence over $A\cup J$. This proves the lemma.
\end{proof}

We now weaken the hypothesis in Lemma~\ref{lem:morley_conj1} from
a Morley sequence to an arbitrary indiscernible sequence.

\begin{lemma}
\label{lem:conj_no_div} Let $|A|<\kappa$ and $p(x,b)$,
$|b|<\kappa$, a type over $b$ that does not divide over $A$. Let
$I$ be an infinite $A-$indiscernible sequence in $tp(b/A)$. Then,
for any $b_0,\,b_1\in I$, $p(x,b_0)\cup p(x,b_1)$ does not divide
over $A$.
\end{lemma}

\begin{proof}
Without loss of generality, $I$ is indexed by $\kappa$. By
Lemma~\ref{lem:all_r_morley}, there is $\alpha<\kappa$ so that
$K=\{a_i:\:\alpha\geq i<\kappa\,\}$ is a Morley sequence over
$A\cup J$, for $J=\{\,a_i:\:i<\alpha\,\}$. For $a\in K$, $p(x,a)$
does not divide over $A$, hence it does not divide over $A\cup J$.
By the Extension Property there is a $c$ realizing $p(x,a)$ such
that $q(x,a)=tp(c/A\cup J\cup\{a\})$ does not divide over $A$.
Since $K$ is a Morley sequence over $A\cup J$, for $a,\,a'\in K$,
$q(x,a)\cup q(x,a')$ does not divide over $A\cup J$. Thus, there
is $d$ realizing $q(x,a)\cup q(x,a')$ such that $tp(d/A\cup
J\cup\{a,a'\})$ does not divide over $A\cup J$. Since $d$ realizes
$q(x,a)$, $tp(d/A\cup J)$ does not divide over $A$. By
Transitivity (Corollary~\ref{cor:trans}), $tp(c/A\cup
J\cup\{a,a'\})$ does not divide over $A$. A fortiori, $p(x,a)\cup
p(x,a')$ does not divide over $A$.
\end{proof}

\begin{lemma}
\label{lem:amalg1} Let $|A|<\kappa$ and $p(x,a)$ and $q(x,b)$ be
types over $A\cup\{a\}$ and $A\cup\{b\}$ respectively. Assume that
$lstp(b/A)=lstp(b'/A)$ and that $tp(a/Abb')$ does not divide over
$A$. If $p(x,a) \cup q(x,b)$ does not divide over $A$, then
$p(x,a) \cup q(x,b')$ does not divide over $A$.
\end{lemma}

\begin{proof}
By Lemma~\ref{l:main}, let $b_0=b$, $b_1, \dots, b_n=b'$ be a
sequence such that there exists $A$-indiscernible sequences
containing $b_i$ and $b_{i+1}$, for $i < n$. By the extension
property, we may find $a'$ realizing $tp(a/Abb')$ such that
$tp(a'/Ab_0 \dots b_n)$ does not divide over $A$. Hence, by using
an automorphism fixing $b$ and $b'$, we may assume that $tp(a/Ab_0
\dots b_n)$ does not divide over $A$. Hence, $tp(a/Ab_ib_{i+1})$
does not divide over $A$, so it is enough to show the conclusion
when $b$, $b'$ belong to the same indiscernible sequence $I$
(since then $p(x,a) \cup q(x,b_1)$ does not divide over $A$, and
so $p(x,a) \cup q(x,b_2)$ does not divide over $A$, etc.).

\begin{claim}
There is an $A$-indiscernible sequence $\{\, a^i b^i :\: i \in
\Bbb{Z}\,\}$ such that $tp(a^i, b^i, b^{i+1}/A) = tp(a, b, b'/A)$,
for each $i \in \Bbb{Z}$.
\end{claim}
Let $\alpha$ be a cardinal, $\pi<\alpha\leq \lambda$, and
$\alpha^*$ the order with $\alpha$ reversed with the elements of
$\alpha^*$ denoted $-i$, for $i\in\alpha$. Since $I$ is
$A$-indiscernible, we may in fact assume that it is of the form
$\{\,b_j:\:j\in\alpha^*+\alpha\,\}$, $b=b_{-0}$ and $b'=b_0$.
Notice that $I'=\{\, b_{-n}\hat{\;}b_n :\: n < \alpha\,\} $ is
also indiscernible over $A$. Since $tp(a/bb'A)$ does not divide
over $A$,
 Lemma~\ref{lem:div_prop1} says we can choose $I$ so that $I'$  is
indiscernible over $A\cup\{a\}$. Notice that
$tp(ab_{-n}/A)=tp(ab/A)$ and $tp(ab_{n}/A)=tp(ab'/A)$, for each
$n<\alpha$. By the indiscernibility of $I$ and the homogeneity of
$M$, for $j\in\alpha^*+\alpha$ there is an $a_j$ such that
$tp(a_jb_k/A)=tp(ab/A)$ for $k\leq j$ and $tp(a_jb_l/A)=tp(ab'/A)$
for $l>j$. By Lemma~\ref{lem:exist_indisc}, there is an
$A$-indiscernible sequence $\{\, a^i b^i :\: i \in \Bbb{Z}\,\}$
such that
$tp(a^0b^0,\dots,a^nb^n/A)=tp(a_{i_0}b_{i_0},\dots,a_{i_n}b_{i_n}/A)$
for some $i_0<\dots<i_n\in\alpha^*+\alpha$. In particular,
$tp(a^ib^ib^{i+1}/A)=tp(abb'/A)$, completing the proof of the
claim.

Since $p(x,a^0)\cup q(x,b^0)$ does not divide over $A$,
$p(x,a^0)\cup q(x,b^0)\cup p(x,a^1)\cup q(x,b^1)$ does not divide
over $A$, by Lemma~\ref{lem:conj_no_div}. In particular,
$p(x,a^0)\cup q(x,b^1)$ does not divide over $A$. Since,
$tp(a^0b^1/A)=tp(ab'/A)$ the lemma is proved.
\end{proof}

\noindent \emph{Proof of Theorem~\ref{thm:tp_amal}: } We seek an
$a$ realizing $lstp(a_i/cb_i)$, for $i=1,2$, such that
$tp(a/cb_1b_2)$ does not divide over $c$. We \emph{claim} that it
is enough to find an $a$ realizing $tp(a_i/cb_i)$, for $i=1,2$,
such that $tp(a/cb_1b_2)$ does not divide over $c$. As we'll see
below, the hypotheses imply that $tp(a_i/cb_i)$ is large. For
$i=1,2$ let $e_i,f_i$ be such that $\{a_i,e_i, f_i\}$ is contained
in a Morley sequence in $tp(a_i/cb_i)$ over $c$. We can choose
these new elements so that $tp(e_1f_1/cb_1b_2e_2f_2)$ does not
divide over $c$. Let $b_i'=b_ie_i$ and consider
$q_i=tp(f_i/cb_i')$. Noting also that $lstp(f_1/c)=lstp(f_2/c)$,
the special case assumed in this claim implies there is an $a$
realizing $q_1\cup q_2$ such that $tp(a/cb_1'b_2')$ does not
divide over $c$. Since $e_1,\,f_1,\,a_1$ are in an indiscernible
sequence over $cb_1$, $lstp(a/cb_1)=lstp(e_1/cb_1)=
lstp(a_1/cb_1)$, and correspondingly over $cb_2$. This proves the
claim.

We first have to deal with the cases when some element is in the
small closure of another. Let $p_i(x,y)=tp(a_ib_i/c)$, for
$i=1,2$. If $a_i\in s\cl(c)$, then $a_1=a_2$ (by
Lemma~\ref{lem:small_lstp}) and we are done. So, we can assume
that $tp(a_1/c)$ is large, hence $tp(a_1/cb_1)$ and $tp(a_2/cb_2)$
are large (Proposition~\ref{prop:main_get_small}). Let $b_2'$ be
such that $lstp(a_1b_2'/c)=lstp(a_2b_2/c)$
(Lemma~\ref{lem:fatten_lstp}). Suppose $b_2'\in s\cl(a_1c)$. Since
$tp(b_2'/a_1c)$ does not divide over $c$, $b_2'\in s\cl(c)$ ,
implying that $b_2'=b_2$, from which the lemma follows easily. So,
we are left with the case when $tp(b_2'/ca_1)$ is large. By
Strong Extension (Lemma~\ref{lem:str_ext}) we can require that
$tp(b_2'/ca_1b_1b_2)$ does not divide over $ca_1$. By
Transitivity, $tp(b_2'/ca_1b_1b_2)$ does not divide over $c$,
hence $tp(b_2'/cb_1b_2)$ does not divide over $c$. By Symmetry and
Transitivity, $tp(b_1/cb_2b_2')$ does not divide over $c$. Also,
by several applications of symmetry and transitivity,
$tp(a_1/cb_1b_2')$ does not divide over $c$. Thus, $p_1(x,b_1)\cup
p_2(x,b_2')$ does not divide over $c$. By Lemma~\ref{lem:amalg1},
$p_1(x,b_1)\cup p_2(x,b_2)$ does not divide over $c$. $\eop$
\bigskip

The most important application of Type Amalgamation is the
connection of parallelism with freeness found in the next section.

\section{Parallelism, imaginary elements and canonical bases}
\label{sec:parallel}

One of the central concepts of geometrical stability theory (for
first-order stable theories) is the notion of a definable family
of  uniformly definable sets. In algebraic geometry a uniform
family of plane curves is a family $\{\,X_{\dd}:\:\dd\in Y\}$ of
one dimensional subsets $X_{\dd}$ of $K^2$, where $X_{\dd}$ is,
 for some polynomial
$f(x,y,\zz)$, the solution set of $f(x,y,\dd)=0$, and $\dd$ ranges
over the elements of the variety $Y\subset K^n$. Of course, for a
smooth theory distinct elements of the family should be almost
disjoint; i.e., have finite intersection. In this case, the
dimension of $Y$ is also called the dimension of the family, a
number which affects the intersection theory of the family and
other properties. This is one reason for normally restricting to
irreducible varieties $X_{\dd}$. Also we require $f$ to be
normalized so that $\dd$ is uniquely determined by $X_{\dd}$.

In a simple model we will have notions of dimension that mimic the
dimension of varieties without a great deal of work. However, new
concepts need to be developed to play the part of the Zariski
topology and the field of definition of a variety. The first
notion is parallelism, an equivalence relation that expresses when
types have a common extension. A canonical base is an equivalence
class of the parallelism relation. With hyperimaginary elements we
have a method of extending the concepts of types and dividing to
such classes. This is required to introduce dimension as a
property of a definable family of definable relations.

Throughout the section $(M,\cR)$ is an $s\lambda-$homogeneous
$\kappa-$simple logical structure.

\begin{definition}
\label{defn:amal_base} A large type $p\in S(A)$ is an
\emph{amalgamation base} if type amalgamation holds for $p$; i.e.,
for any $B,\,C$, with $|B|,\,|C|<\kappa$ and $B$ free from $C$
over $A$, and types $q$ over $B$ and $r$ over $C$ such that $p$
has a free extension over $A\cup B$ containing $q$ and a free
extension over $A\cup C$ containing $r$, $p$ has a free extension
over $A\cup B\cup C$ containing $q\cup r$.
\end{definition}

\begin{remark}
\label{rem:defn_amal_base} When $|A|<\kappa$, a large $p\in S(A)$
is an amalgamation base if and only if for any $a,\,b$ realizing
$p$, $lstp(a/A)=lstp(b/A)$. [The direction from right to left is
virtually a restatement of the Type Amalgamation Theorem. To prove
the other direction, suppose that $a$ and $b$ are any two
realizations of $p$ such that $tp(a/bA)$ does not divide over $A$.
There are types $q(x,a)$ and $r(x,b)$ over $A\cup\{a\}$ and
$A\cup\{b\}$, respectively, which do not divide over $A$, such
that if $c$ realizes $q(x,a)$ and $d$ realizes $r(x,b)$, then
$lstp(c/A)=lstp(a/A)$ and $lstp(d/A)=lstp(b/A)$. Since $p$ is an
amalgamation base there is a $c$ realizing $q(x,a)\cup r(x,b)$.
Thus, $lstp(a/A)=lstp(b/A)$. ]
\end{remark}

\subsection{Parallelism}
\label{subsec:paral}

In a stable theory stationary types are parallel if they have a
common free extension. In a simple homogeneous model (or even a
simple theory) the role of parallelism is played by a slightly
more complicated concept about amalgamation bases.

\begin{definition}
\label{defn:parallel} Given $p\in S(A)$ and $q\in S(B)$
amalgamation bases  we write $p\sim_1 q$ if $p\cup q$ does not
divide over $A$ and does not divide over $B$. We write $p\sim q$
and say $p$ is \emph{parallel} to $q$ if there are amalgamation
bases $q_0,\dots,q_k$ such that $p=q_0\sim_1 q_1\sim_1\dots\sim_1
q_k=q$.
\end{definition}

\begin{remark}
\label{rem:defn_parallel} Parallelism is clearly an equivalence
relation on the class of amalgamation bases. In fact it is the
transitive closure of $\sim_1$. Parallelism is invariant under
automorphisms.
\end{remark}

Suppose $p\in S(a)$ is an amalgamation base and $p'\in S(a')$ is
conjugate to $p$, where $a$ is free from $a'$ over $\emptyset$. If
$p$ does not divide over $\emptyset$ then Type Amalgamation says
that $p\cup p'$ does not divide over $\emptyset$ (assuming there
are $b$ realizing $p$ and $b'$ realizing $p'$ with
$lstp(b)=lstp(b')$). In this case $p\sim p'$. This suggests that
``most'' conjugates are parallel if $p$ does not divide
 over $\emptyset$.
This is formalized in the next result, which guides all  uses of
parallelism.

\begin{proposition}
\label{prop:par_free} Let $p\in S(a)$ be an amalgamation base with
$|a|<\kappa$, $b$ a sequence
 and $\mathbf{P}$ the class of types conjugate to
$p$ over $b$. Then $p$ does not divide over $b$ if and only if
$\mathbf{P}$ contains a bounded number of parallelism classes.
\end{proposition}

\begin{proof}[\emph{of the left to right direction}]
Suppose that $p$ does not divide over $b$. There is a subsequence
$b_0$ of $b$ of length $<\kappa$ such that $p$ does not divide
over $b_0$ so we may as well assume $|b|<\kappa$. For $q$ a Lascar
strong type over $b$ let $\mathbf{P}_q$ be the set of
$p'\in\mathbf{P}$,  $p'\in S(a')$ such that there is $c$ realizing
$p'\cup q$ and $tp(c/a'b)$ does not divide over $b$.
\begin{claim}
If $p_0,\,p_1\in\mathbf{P}_q$, for some Lascar strong type $q$
over $b$, then $p_0\sim p_1$.
\end{claim}
Given $p_i\in S(a_i)$ in $\mathbf{P}_q$, for $i=0,1$, let $r\in
S(c)$ in $\mathbf{P}_q$ such that $tp(c/ba_0a_1)$ does not divide
over $b$. By Type Amalgamation there are, for $i=0,1$, $c_i$
realizing $p_i\cup r$ such that $tp(c_i/cba_i)$ does not divide
over $b$. Thus, $p_i\sim_1 r$, for $i=0,1$, hence $p_0\sim p_1$,
proving the claim.

Since the class of Lascar strong types over $b$ is bounded,
$\mathbf{P}$ contains a bounded number of parallelism classes.
\end{proof}

The opposite direction of the proposition requires some detailed
analysis of the behavior of amalgamation bases. Recall that if
$(M,\cR)$ is almost $\kappa-$simple, $\indep{A}{C}{B}$ denotes $A$
is free from $B$ over $C$.

\begin{lemma}
\label{lem:amalg_base_chng} Let $r\in S(d)$ be an amalgamation
base, where $|d|<\kappa$, $I$ a  Morley sequence in $r$ over $d$
indexed by $\kappa^*$ and $c$ a realization of $r$ free from $I$
over $d$. Then, $c$ is free from $d$ over $I$.
\end{lemma}

\begin{proof}
In the special case where $I\hat{\;}\langle c\rangle$ is
$d-$indiscernible, this follows basically from
Lemma~\ref{lem:all_r_morley}. In general we will use this fact and
that $r$ is an amalgamation base. Recall Remark~\ref{rem:indep}
about independence of Morley sequences.
\begin{claim}
There is a sequence $J$ of length $<\kappa$ such that $I\hat{\;}J$
is $d-$indiscernible, $I$ is independent over $J$ and $d$ is free
from $I\hat{\;}J$ over $J$.
\end{claim}
To see this, first let $I'$ be a sequence indexed by $\kappa^*$
such that $I\hat{\;}I'$ is indiscernible over $d$. Consider $I'$
as an indiscernible sequence in the reverse order (indexed by
$\kappa$) and apply Lemma~\ref{lem:all_r_morley} to obtain a
sequence $J\subset I'$ of length $<\kappa$ such that, letting
$i_0=\inf\{\,j:\:a_j\in J\,\}$, $I''=\{\,a_i\in I':\:i<i_0\,\}$ is
a Morley sequence over $J$ (under the reverse order). Thus, $I''$
is independent over $J$. Since $I\hat{\;}I''$ is indiscernible
over $J$, $I$ is independent over $J$. Now suppose $\aa$ is a
finite subset  of $\bigcup I\cup\bigcup J$ and suppose, towards a
contradiction, that $tp(\aa/J\cup d)$ divides over $J$. There is a
$K\subset I$, $|K|<\kappa$, such that $d$ is free from $I\cup J$
over $K\cup J$. By the indiscernibility of $I\hat{\;}J$ over $d$
we can assume $\aa$ is disjoint from $K$, hence free from $K$ over
$J$. Since $\aa$ is free from $d$ over $J\cup K$ transitivity of
independence implies that $\aa$ is free from $d$ over $J$. This
contradiction proves the claim.

Arguing as in the claim there is also a (nonempty) $K\subset I$,
$|K|<\kappa$, such that $I\cup J$ is free from $d$ over $K.$ Let
$s(x,z)=tp(J,K/d)$. Let $c'$ be any element of $K$ and
$q(x,y)=tp(J,c'/d)$.

\begin{claim}
There is $J''$ free from $K\cup c$ over $d$ such that $q(J'',c)$
and $s(J'',K)$ hold.
\end{claim}
We first show there is $J'$ realizing $lstp(J/d)$ and $q(x,c)$.
Since $r$ is an amalgamation base, $lstp(c/d)=lstp(c'/d)$. So,
there is $f$, a strong automorphism over $d$ taking $c'$ to $c$.
Then $lstp(f(J)/d)=lstp(J/d)$ and $q(f(J),c)$ holds. Since $J$ is
free from $c'$ over $d$ and $K$ is free from $c$ over $b$ Type
Amalgamation yields $J''$ free from $K\cup c$ over $d$ such that
$q(J'',c)$ and $s(J'',K)$ hold.

The following chain of arguments shows that $\indep{d}{K}{c}$.
$\indep{J''}{d}{K\cup\{c\}}\;\Longrightarrow\;
\indep{J''}{K\cup\{d\}}{c}$ (transitivity)
$\Longrightarrow\;\indep{c}{K\cup\{d\}}{J''}$ (symmetry)
$\Longrightarrow\;\indep{c}{d}{J''\cup K}$ (transitivity and the
independence of $c$ and $K$ over $d$)
$\Longrightarrow\;\indep{c}{J''\cup \{d\}}{K}$ (transitivity)
$\Longrightarrow\;\indep{c}{J''}{K\cup \{d\}}$ (transitivity and
the fact that $q(J'',c)$ holds) $\Longrightarrow\;\indep{c}{K\cup
J''}{d}$
 (transitivity) $\Longrightarrow\;\indep{d}{K\cup J''}{c}$
 (symmetry) $\Longrightarrow\;\indep{d}{K}{\{c\}\cup J''}$
 (transitivity and the fact that $s(J'',K)$ holds)
 $\Longrightarrow\;\indep{d}{K}{c}$ (a fortiori).

Finally, to show that $\indep{d}{I}{c}$, it suffices to prove that
$\indep{d}{K}{I\cup\{c\}}$ by transitivity and the independence of
$d $ and $I$ over $K$. It is enough to show that
$\indep{d}{K}{I'\cup\{c\}}$, for any $I'\subset I$ with
$|I'|<\kappa$. The argument above can be repeated with $K\cup I'$
replacing $K$ to prove that $\indep{d}{K\cup I'}{c}$. An
application of transitivity then completes the proof.
\end{proof}

\begin{lemma}
\label{lem:par_imp_base} If $q\in S(c)$ is an amalgamation base,
$q'\in S(c')$ is conjugate to $q$ and $q\sim q'$, then there is
$a$ realizing $q'$ such that $\indep{a}{c'}{c}$ and
$\indep{a}{c}{c'}$.
\end{lemma}

\begin{proof}
By the definition of parallelism there  are $q_i$, $i\leq k$, such
that $q=q_0\sim_1\dots\sim_1 q_{k-1}\sim_1 q_k=q'$. The proof is
by induction on $k$, so we assume there is $r\in S(b)$ conjugate
to $q$, $r\sim_1 q'$, and there is $d$ realizing $r$ such that
$tp(d/cb)$ does not divide over $b$ and does not divide over $c$.
Since $r\sim_1 q'$ there is $a$ realizing $r\cup q'$ such that
$\indep{a}{c'}{b}$ and $\indep{a}{b}{c'}$. Without loss of
generality, $\indep{a}{b}{cc'}$. To prove that $\indep{a}{c}{c'}$
a Morley sequence in $r$ must be introduced and the preceding
lemma applied.

Since $r$ does not divide over $c$ there is a sequence $I$ of
length $\kappa$ which is a Morley sequence in $r$, is free from
$c$ over $b$, and free from $b$ over $c$. Without loss of
generality, $I$ is free from $ac'$ over $bc$. Thus, $a$ is free
from $Icc'$ over $b$ (by symmetry and transitivity), and $a$ is
free from $cc'$ over $Ib$. $a$ is free from $b$ over $I$ by
Lemma~\ref{lem:amalg_base_chng}. Also by transitivity, $a$ is free
from $cc'b$ over $I$, hence $a$ is free from $c'b$ over $Ic$.
Since $I$ is free from $a$ over $c$, $a$ is free from $I$ over
$c$, hence $a$ is free from $Ibc'$ over $c$. In particular, $a$ is
free from $c'$ over $c$, proving the lemma.
\end{proof}

\noindent \emph{Proof of Proposition~\ref{prop:par_free} [of right
to left direction]. } Since there are boundedly many parallelism
classes in $\mathbf{P}$ it contains $p'\in S(a')$ and $p''\in
S(a'')$ such that $p'\sim p''$ and $\indep{a'}{b}{a''}$. By
Lemma~\ref{lem:par_imp_base} there is $c$ realizing $p''$ such
that $\indep{c}{a''}{a'}$ and $\indep{c}{a'}{a''}$. Taking a free
extension of $tp(c/a'a'')$ if necessary we can require that $c$ is
free from $b$ over $a'a''$. Since $c$ is free from $a''$ over
$a'b$, $a''$ is free from $c$ over $a'b$. Using
$\indep{a'}{b}{a''}$ and transitivity, $a''$ is free from $c$ over
$b$, hence $c$ is free from $a''$ over $b$. Since $tp(c/a''b)$ is
a free extension of $p''$, $p''$ does not divide over $b$,
completing the proof. Since $p''$ is conjugate to $p$ over $b$,
$p$ does not divide over $b$. $\eop$
\bigskip

After defining imaginary elements the notion of a canonical base
of an amalgamation base will be defined as the parallelism class
of an amalgamation base. Expanding $M$ by the addition of
imaginary elements allows us to treat canonical bases as elements
of the expanded model. Virtually every result in geometrical
stability theory involves canonical bases.

\subsection{Imaginary and hyperimaginary elements}
\label{subsec:imags}

Imaginary elements were introduced by Shelah to enable us to work
with the classes of first-order definable equivalence relations as
if they were elements of the model. Here, equivalence relations
that may not be first-order definable are needed, so some
discussion of the concept is warranted. Hyperimaginary elements
were introduced in \cite{hakipi} to capture canonical bases in all
simple theories.  In that setting a hyperimaginary element is a
class of a type-definable equivalence relation in possibly
infinitely many variables.

In general there is no reason to think the canonical base exists
as a tuple from the model, as it does in an algebraically closed
field. In a stable theory this deficiency is removed by expanding
from $M$ to $M^{eq}$, which contains the classes of any definable
equivalence relation in $M$. When $M$ is saturated so is $M^{eq}$,
thus the expanded universe satisfies the conditions under which
stability theory is developed. In fact, it is standard to simply
work in an expanded universe containing classes for each definable
equivalence relation. In a simple theory parallelism may only be
type-definable. When the classes of a type-definable equivalence
relation are added to a saturated model the resulting model may
not be saturated, hence the classical theory cannot be applied
here. This is handled with so-called hyperimaginary elements
\cite{hakipi}.

In a simple homogeneous model that isn't saturated it isn't clear
that we can even add classes for equivalence relations in $\cR$
and obtain a homogeneous model. Even in the settings where this is
possible, the most important equivalence relation, parallelism,
may not be $\cR-$type definable much less an element of $\cR$.
Thus, the only alternative is to handle all invariant equivalence
relations as we do type-definable equivalence relations in the
first-order case by hyperimaginaries. Admittedly, there is
probably little to be gained by considering classes of invariant
equivalence relations rather than the tuples in the original
model. However, we feel it is worth stating the definitions to
create a framework for more useful results under additional
assumptions. It is still worthwhile to  separate the finitary
equivalence relations from the infinitary ones.

\begin{definition}
\label{defn:imags} Let $(M,\cR)$ be a $s\lambda-$homogeneous
logical structure. Let $\mathcal{E}$ be the class of relations $E$
on $M$ such that (1) $E$ is an equivalence relation on $M^n$, for
some $n<\omega$, and (2) $E$ is invariant under automorphisms of
$(M,\cR)$. For each $E\in\mathcal{E}$ let $S_E$ be a new sort and
let $M_E=M^n/E$ (if $E$ is an equivalence relation on $M^n$) and
$f_E$ the quotient map from $M^n$ into $M_E$. Given $L$ the
language of $M$, $L^{eq}$ is the language obtained from $L$ by
adding the sorts $S_E$ and symbols for the maps $f_E$, for each
$E\in\mathcal{E}$. Let $M^{eq}$ be the expansion of $M$ obtained
by letting $M_E$ interpret the sort $S_E$ and the relevant
quotient map interpret $f_E$. The elements of $M^{eq}$ are called
\emph{imaginary elements}.

It remains to equip  $M^{eq}$ with the relations under which  it
will be a logical structure with types equivalent to orbits for
finite tuples. For each $p=tp(\aa)$, $\aa$ a finite tuple from
$M$, let $R_p=\{\,\bb:\:p(\bb)\,\}$. Let ${\cR}^c$, the
\emph{completion of $\cR$}, be $\cR\cup\{\,R_p:\:p=tp(\aa)\,\}$.
Of course, $(M,\cR)$ and  $(M,{\cR}^c)$ have the same types. On
$M^{eq}$ the structure induced by $(M,\cR)$ is defined as follows.
Let $\mathcal{S}$ be the class of all relations $\bar{R}$, for
$R\in{\cR}^c\cup\mathcal{E}$ defined as follows. Suppose
$E_i(\yy_i,\zz_i)\in\mathcal{E}$, for $i=1,\dots,n$, and
$R(\xx_1,\dots,\xx_k)\in{\cR}^c$ with $\xx_i\subset \yy_i$,
$i=1,\dots,n$. Let $\bar{R}$ be the relation that holds on
$a_1,\dots,a_k$, where $a_i\in M_{E_i}$, if there are
$c_1,\dots,c_k$ such that $c_i/E_i=a_i$ and $R(c_1,\dots,c_k)$
holds (after assigning the corresponding subsequence of $c_i$ to
$\xx_i$). Finally, let ${\cR}^{eq}$ be the closure of
${\cR}^c\cup\mathcal{E}\cup\{\,f_E:\:E\in\mathcal{E}\,\}\cup\mathcal{S}$
under finite unions and intersections.

Fix $\pi$ the cardinal defined in $(\bold{P})$ at the beginning of
Section~\ref{sec:div}. Let $\mathcal{E}^{h}$ be the collection of
equivalence relations $E$ on $M^\mu$, for some $\mu\leq\pi$, such
that $E$ is invariant under automorphisms of $(M,\cR)$. For $E\in
\mathcal{E}^{h}$ let $M_E$ be the set of equivalence classes,
$M^\mu/E$, where $E$ is on $M^\mu$, and $f_E$ the quotient map
from $M^\mu$ into $M_E$. Let $M^h$ be obtained from $M$ be adding
$M_E$ and $f_E$ for each $E\in \mathcal{E}^{h}$. Let $\mathcal{S}$
be the class of all relations $\bar{R}$, for
$R\in\cR\cup{\mathcal{E}}^h$, defined as for $M^{eq}$. Let
${\cR}^h$ be the closure of
$\cR\cup{\mathcal{E}}^h\cup\mathcal{S}$ under finite unions and
intersections. The structure $(M^h,{\cR}^h)$ is the
\emph{expansion of $M$ by hyperimaginaries}.

For $A\cup\{a\}\subset M^h$, $tp(a/A)$ is $\{\,S(x,\bb):\:
S\in\mathcal{S},\;\bb\subset A,\;\textrm{ and }S(a,\bb)\textrm{
holds}\,\}$. In $(M^{eq},{\cR}^{eq})$, $tp(a/A)$ is simply the
${\cR}^{eq}-$type in the logical structure $(M^{eq},{\cR}^{eq})$
(see Lemma~\ref{lem:first_meq}).
\end{definition}

\begin{remark}
\label{rem:defn_meq} (i) In the first-order definition of
$M^{eq}$, classes are only added for first-order definable
equivalence relations. Classes of type-definable equivalence
relations in finitely many variables are normally called finitary
hyperimaginaries (see \cite{hakipi}). In this setting where
compactness isn't used all invariant equivalence relations in
finitely many variables can be treated alike. Of course, in some
applications the distinction between type-definable and invariant
equivalence relations may be important in analyzing particular
structures. Invariant equivalence relations in infinitely many
variables are very problematic without compactness, a point we
will expand on below.

(ii) Since $E\in\mathcal{E}$ is invariant under automorphisms of
$M$, $\{\,ab:\:E(a,b)\,\}$ is a union of orbits of finite tuples.
So, if $(M,\cR)$ has $\mu$ complete $\cR-$types,
$|\mathcal{E}|\leq 2^\mu$. This bounds the number of sorts added
to $L$ to form $L^{eq}$.
\end{remark}

\begin{lemma}
\label{lem:first_meq} (i) Any automorphism of $M$ extends to an
automorphism of $M^{eq}$ and each automorphism of $M^{eq}$ is
determined by its restriction to $M$.

(ii) $(M^{eq},\cR^{eq})$ is a logical structure.

(iii) Given $a,\,a',\,\,b\in M^{eq}$, $tp(a/b)=tp(a'/b)$ if and
only if there is an automorphism $f$ of $M^{eq}$ fixing $b$ with
$f(a)=a'$.
\end{lemma}

\begin{proof}
(i) Since any element of $M^{eq}$ is the image of a tuple from $M$
under a function interpreting a symbol in $L^{eq}$, this is
straightforward.

(ii) By (i) the elements of $\cR^{eq}$ are invariant under
automorphisms of $M^{eq}$. The other conditions follow quickly
from the construction of $\cR^{eq}$ and the fact that $(M,\cR)$ is
a logical structure.

(iii)  By the construction of $\cR^{eq}$, for any finite tuple $a$
from $M$ there is $R\in\cR$ such that $R(b)$ holds if and only if
$tp(b)=tp(a)$. For the given elements $a,\,a',\,b$ suppose
$a,\,a'\in M_E$ and $b\in M_F$. Let $a=c/E$, $b=d/F$ and
$S(xy)\in\cR$ define $tp(cd)$. The relation $\bar{S}(uv)$
expressing $\exists xy(\;x/E=u\wedge y/F=v\wedge S(xy)\;)$ is in
$tp(ab)$, hence in $tp(a'b)$. Thus, there are $c',\,d'\in M$,
$c'/E=a'$, $d'/F=b$ and $S(c'd')$. There is an automorphism $f$ of
$M$ (and $M^{eq}$) such that $f(cd)=c'd'$, hence $f(b)=b$ and
$f(a)=a'$, proving the lemma.
\end{proof}

\begin{notation}
For $a\in M^h$ and $A\subset M^h$, we write $a\in dcl(A)$ if $a$
is fixed by every automorphism of $(M^h,{\cR}^h)$ that fixes $A$
pointwise. We say $a$ is in the \emph{bounded closure of $A$},
$a\in bdd(A)$, if the orbit of $a$ under the group of
automorphisms that fixes $A$ pointwise is bounded (i.e, of
cardinality $<\lambda$). The \emph{length of $a$}, denoted $|a|$,
is the least cardinal $\nu$ such that $a\in dcl(B)$ for some
$B\subset M$ with $|B|=\nu$.
\end{notation}

\begin{definition}
Given a sequence $\aa$ from $M^{eq}$ of length $\alpha$ for some
ordinal $\alpha<|M|$, a \emph{base} for $\aa$ is a sequence
$\{\,c_i:\:i<\alpha\,\}\subset M$ such that, if $a_i$ is in
$M_{E_i}$, $a_i=c_i/E_i$.
\end{definition}

It is an open problem to determine when $(M^{eq},\cR^{eq})$ is
$s\lambda-$homogeneous. Here is the problem.   The
$s\lambda-$homogeneity of $(M^{eq},\cR^{eq})$ boils down to
showing that for sequences $\aa,\,\bb$ from $M^{eq}$ of length
$<\lambda$, if $tp(\aa)=tp(\bb)$ then there are bases $\cc,\,\dd$
of $\aa$, respectively, such that $tp(\cc)=tp(\dd)$. This can be
proved for type-definable equivalence relations when $M$ is
compact, but it is generally open. Similarly, for $(M^h,{\cR}^h)$.
We state for the record:

\begin{lemma}
\label{lem:compact_imag} Suppose $(M,\cR)$ is a compact
$s\lambda-$homogeneous logical structure (i.e., a large saturated
model of its complete theory). Let $E$ be an equivalence relation
that is equivalent to an $\cR-$type on $(M,\cR)$. Let
$(M',{\cR}')$ be the restriction of $(M^{eq},{\cR}^{eq})$ to
$M\cup M_E$. Then $(M',{\cR}')$ is $s\lambda-$homogeneous.
\end{lemma}

\begin{proof}
Left to the reader with the hint in the preceding paragraph.
\end{proof}

Even if $(M^{eq},\cR^{eq})$ is not $s\lambda-$homogeneous we can
discuss dividing in it. The treatment is similar to the handling
of hyperimaginaries in \cite{hakipi}. For example, indiscernible
sequences of imaginary elements exist as they do for ordinary
elements in the $s\lambda-$homogeneous model $(M,\cR)$. In fact,
the treatment for $(M^{eq},\cR^{eq})$ and $(M^h,{\cR}^h)$ is the
same so we work in the more general structure.

\begin{lemma}
\label{lem:imag_indisc} Let $A\subset M^h$, $a\in M^h$ such that
$tp(a/A)$ is large and $X$ any order type of length $\leq
\lambda$. Then $M^h$ contains an $A-$indiscernible sequence
$I=\{\,a_i:\:i\in X\,\}$ in $tp(a/A)$.
\end{lemma}

\begin{proof}
Given $E$ so that $a\in M_E$ there is a $b\in M$ such that $b/E=a$
and $q=tp(b/A)$ is large. Let $J=\{\,b_i:\:i\in X\,\}$ be an
$A-$indiscernible sequence in $q$, whose existence is guaranteed
by the $s\lambda-$homogeneity of $M$. Then, $I=\{\,b_i/E:\:i\in
X\,\}$ is an $A-$indiscernible sequence in $tp(a/A)$.
\end{proof}

A little extra care is required when dealing with indiscernibles
of imaginaries. Most of the subtleties are revealed in the
following definition and lemma.

\begin{definition}
An infinite $A-$indiscernible sequence $I=\{\,a_i:\:i\in X\,\}$ is
\emph{ultra-indiscernible over $A$} if there is a base $A'$ for
$A$ and an $A'-$indiscernible sequence $J=\{\,c_i:\:i\in X\,\}$
such that $a_i=c_i/E$.
\end{definition}

\begin{lemma}
\label{lem:imag_indisc2} Let $A$ be a set and $I=\{\,a_i:\:i\in
X\,\}$ an infinite $A-$indiscernible sequence in  $(M^h,{\cR}^h)$.
If $X$ is sufficiently large then there is an infinite $J$ which
is ultra-indiscernible over $A$ and has the same type diagram over
$A$ as $I$. If $I=\{\,c_i/E:\:i\in X\,\}$ and $I'=\{\,c_i:\:i\in
X\,\}$ are $A-$indiscernible, then for any ordered set $Y$,
$|Y|<|M|$, there are $a_i$, for $i\in Y$, such that
$J=\{\,a_i:\:i\in X+Y\,\}$ is $A-$indiscernible.
\end{lemma}

\begin{proof}
Let $A'$ be a base for $A$. Suppose $a_i\in M_E$ and let
$a_i=c_i/E$. If $X$ is sufficiently large there is an infinite
$A'-$indiscernible sequence $J_0=\{\,d_j:\:j<\omega\,\}$ such that
for all $n<\omega$, there are $i_1<\dots<i_n\in X$ with
$tp(d_1,\dots,d_n/A')= tp(c_{i_1},\dots,c_{i_n}/A')$.
>From here the proof is as in Lemma~\ref{lem:large_indisc}.
\end{proof}

Dividing is defined as it is in $(M,\cR)$ with an additional
condition on the indiscernibles.

\begin{definition}
\label{defn:divide_imag} Let $p(v,b)$ be a type over the element
$b$ in $(M^h,{\cR}^h)$, and $A$ a set. Then $p(v,b)$ \emph{divides
over $A$} if there is a sequence $\{\,b_i:\:i\in X\,\}$,
$|X|<|M|$,  ultra-indiscernible over $A$ with $tp(b_i/A)=tp(b/A)$,
such that $\bigcup_{i\in X}p(v,b_i)$ is inconsistent.
\end{definition}

The connection between dividing in $M$ and dividing in $M^h$ is
the following.

\begin{lemma}
\label{lem:connect_div} Let $p(v,b)=tp(a/b)$ be a type over the
element $b$ in $(M^h,{\cR}^h)$, and $A$ a set. For $a\in M_E$ and
$b\in M_F$ the following are equivalent.
\begin{enumerate}[(1)]
\item $p(v,b)$ divides over $A$. \item There is a base $A'$ for
$A$, such that for all $c$ and $d$ with $c/E=a$ and $d/F=b$,
$tp(c/d)$ divides over $A'$.
\end{enumerate}
\end{lemma}

\begin{proof}
$(1)\Longrightarrow(2)$ First observe that $tp(b/A)$ is large,
since $p(v,b)$ divides over $A$. Let $J=\{\,e_i/F:\:i\in X\,\}$ be
an infinite sequence ultra-indiscernible over $A$ in $tp(b/A)$,
where $\{\,e_i:\:i\in X\,\}$ is indiscernible over the base $A'$
for $A$ and $\bigcup_{i\in X}p(v,e_i/F)$ is inconsistent in
$(M^h,{\cR}^h)$. To reach a contradiction, suppose there are
$c,\,d$ such that $c/E=a$, $d/F=b$ and $q(v,d)=tp(c/d)$ does not
divide over $A'$. Since $d/F=b$, for each $i\in X$ there is $f_i$
realizing $tp(d/A')$ such that $f_i/F=e_i/F$. By standard
arguments about indiscernibles there are $f_i$, $i\in X$, such
that $K=\{\,f_i:\:i\in X\,\}$ is $A'-$indiscernible in $tp(d/A')$
and $f_i/F=e_i/F$. Since $tp(c/d)$ does not divide over $A'$ there
is $c'$ realizing $\bigcup_{i\in X}q(v,f_i)$. Then, $c'/E$
realizes $\bigcup_{i\in X}p(v,e_i/F)$. This contradiction
completes the proof of this direction.

$(2)\Longrightarrow(1)$ Now suppose that $p(v,b)=tp(a/b)$ does not
divide over $A$, $A'$ is any base for $A$ and $d/F=b$. We need to
find a $c$ such that $c/E=a$ and $tp(c/d)$ does not divide over
$A'$. If $tp(d/A')$ is small then for any $c$ such that $c/E=a$,
$tp(c/d)$ does not divide over $A'$. So we can assume $tp(d/A)$ is
large and let $J=\{\,d_i:\:i\in X\,\}$ be a Morley sequence in
$tp(d/A')$ where $X$ satisfies $(\ast)$ in
Proposition~\ref{prop:witness_divide}. Then, $J'=\{\,d_i/F:\:i\in
X\,\}$ is ultra-indiscernible over $A$ in $tp(b/A)$. Since
$p(v,b)$ does not divide over $A$ there is $a'$ realizing
$\bigcup_{i\in X}p(v,d_i/F)$. Let $a'=c'/E$. $X$ can be chosen
arbitrarily large so there is an arbitrarily large $Y\subset X$
also satisfying $(\ast)$ in Proposition~\ref{prop:witness_divide}
such that $tp(c'd_i)=tp(c'd_j)$, for $i,\,j\in Y$. Thus, for
$q(xy)=tp(c'd_i)$, $i\in Y$, Proposition~\ref{prop:witness_divide}
shows that $q(x,d)$ does not divide over $A$. Since there is a $c$
realizing $q(x,d)$ such that $c/E=a$, the proof is complete.
\end{proof}

Suppose $(M,\cR)$ is $\kappa-$simple. With this translation of
dividing in $M^h$ to dividing in $M$ the key properties of
dividing can be shown to extend to $M^h$. In particular, symmetry
and transitivity follow from Lemma~\ref{lem:connect_div} by a
straightforward argument. Furthermore, given $A\subset M^h$ and
$a$ a sequence of length $<\kappa$ in $M^h$, there is a $b\subset
A$ of length $<\kappa$ such that $a$ is free from $A$ over $b$.
Also, for $tp(a/A)$ large  in $M^h$ there are arbitrarily long
Morley sequences in $tp(a/A)$.

For $a,\,b\in M^h$ we say \emph{$a$ and $b$ have the same Lascar
strong type over $A$}, written $lstp(a/A)=lstp(b/A)$, if there are
$a',\,b'\in M$ and $E$ and base $A'$ for $A$ such that $a'/E=a$,
$b'/E=b$ and $lstp(a'/A')=lstp(b'/A')$. If $(M^{eq},\cR^{eq})$ is
$s\lambda-$homogeneous this agrees with the
Definition~\ref{defn:lstp}.
\bigskip

\subsection{Canonical bases}
\label{subsec:can_base}

Here classes of the parallelism relation are captured as
hyperimaginary elements. This ``internalizes'' parallelism and
allows us to study the properties of canonical bases with respect
to dividing.

\begin{definition}
\label{defn:can_base} Let $(M,\cR)$ be $\kappa-$simple and $p\in
S(a)$ an amalgamation base. A \emph{canonical base} for $p$ is a
maximal hyperimaginary $c$ such that for any automorphism
$\alpha$, $\alpha(p)$ is parallel to $p$ if and only if
$\alpha(c)=c$.
\end{definition}

\begin{lemma}
\label{lem:can_base_exist} Let $(M,\cR)$ be $\kappa-$simple and
$p(v,a)\in S(a)$ an amalgamation base with $|a|<\kappa$. There is
a unique canonical base $c$ for $p$ and $c=dcl(c_0)$ for some
hyperimaginary $c$ with $|c|<\kappa$.
\end{lemma}

\begin{proof}
Let $E$ be the equivalence relation on $tp(a)$ such that $E(b',b)$
holds if and only if $p(v,b')$ is parallel to $p(v,b)$. Let $c_0$
be $a/E$. Let $c$ be the unique hyperimaginary element obtained by
adjoining all hyperimaginaries in $dcl(c_0)$. Clearly, $c$ is the
unique canonical base for $p$.
\end{proof}

The most fundamental property of a canonical base is

\begin{proposition}
\label{prop:can_base_div} Let $(M,\cR)$ be $\kappa-$simple,
$p(v,a)\in S(a)$ an amalgamation base, $|a|<\kappa$, and $c$ the
canonical base for $p$. For any set $B$, $p(v,a)$ does not divide
over $B$ if and only if $c\in bdd(B)$.
\end{proposition}

\begin{proof}
This is immediate by Proposition~\ref{prop:par_free}.
\end{proof}

\section{Stability}
\label{sec:stab}

An $s\lambda-$homogeneous model $M$ is \emph{$\mu-$stable} if for
any $A\subset M$ with $|A|\leq\mu$, $M$ realizes $\leq\mu$
complete types over $A$. An $s\lambda-$homogeneous model is {\em
stable} if it is $\mu$-stable for some $\mu$ ($<|M|$). As we
pointed out in the first section, stability of
$s\lambda-$homogeneous models was initiated by Shelah and has been
studied extensively by Shelah, Hyttinen, Grossberg, Lessmann and
others. A first-order theory is stable if and only if it is simple
and there is a uniform bound on the number of free extensions of a
complete type. For $s\lambda-$homogeneous models, this fails:
stability does not imply simplicity. An example due to Shelah of a
stable homogeneous model where the free extension property fails
can be found in \cite{HyLe}; a description of this example is
given in the next section. We will therefore study stability under
the additional assumption of simplicity. We make the following
definition.

\begin{definition}
\label{defn:simply_stable} An  $s\lambda-$homogeneous model
$(M,\cR)$ is \emph{$\kappa$-simply stable in $\mu<|M|$} if
$(M,\cR)$ is $\kappa$-simple and if $p$ is a complete $\cR-$type
over $A\subset M$, $|A|\leq\lambda$ then for any $B$, $A\subset
B\subset M$ and $|B|\leq\mu$, there are $\leq\mu$ complete types
over $B$ which are free extensions of $p$. $(M,\cR)$ is
\emph{$\kappa-$simply stable} if it is $\kappa-$simply stable in
some $\mu<|M|$; and \emph{simply stable} is defined similarly.
\end{definition}

\begin{remark}
\label{rem:defn_simp_stab} If $M$ is the universal domain of a
simple first-order theory and $\cR$ is the class of definable
relations, then $(M,\cR)$ is simply stable if and only if $Th(M)$
is stable; i.e., $(M,\cR)$ is stable as a $s\lambda-$homogeneous
model.
\end{remark}

We recall a few facts about stability of $s\lambda-$homogeneous
models. Suppose $(M,\cR)$ is a stable $s\lambda-$homogeneous
model. We will use the notation of logical structures where
Shelah, \textit{et al}, consider structures in a first-order
language, however, the proofs are the same and although we quote
the general results, we will only use particular cases.

If $(M, \cR)$ is stable, then there is a first cardinal, written
$\lambda(M)$, for which $(M,\cR)$ is stable in $\lambda(M)$. A
complete type $tp(a/A)$ is said to {\em split strongly over $B$},
if there is an $R\in\cR$ and an infinite $B$-indiscernible
sequence $\{\,b_i:\:i < \omega\,\}$ such that $b_0, b_1 \in A$,
then $R(x,b_0) \in p$ if and only if $R(x,b_1) \in p$. There is a
least cardinal, written $\kappa(M)$, such that for each finite
$a$, and set $C$, there exists $B \subset C$ of size less than
$\kappa(M)$ such that $tp(a/C)$ does not split strongly over $B$.
In \cite{Sh:4} (although with these definitions this is done in
\cite{GrLe}), Shelah proves that, with $T=Th(M)$, $\kappa(M) \leq
\lambda(M) < \beth_{(2^{|T|})^+}$ and that $M$ is stable in $\mu$
if and only if $\mu \geq \lambda(M)$ and $\mu^{< \kappa(M)} =
\mu$. This is the {\em Stability Spectrum Theorem}. Also in
\cite{Sh:4}, Shelah showed that if $M$ is stable in $\mu$, $I$ is
a set of size $\mu^+$ containing finite sequences, and $A$ has
size at most $\mu$, then there is an $A$-indiscernible subset of
$I$ of size $\mu^+$.

If $I$ is an $A$-indiscernible sequence (hence set), with $\ell(a)
< \kappa(M)$ for each $a \in I$, then if $b$ has size less than
$\kappa(M)$, there exists $J \subset I$ of size less than
$\kappa(M)$ such that $I \setminus J$ is $Ab$-indiscernible (see
\cite{GrLe}). In particular, for an indiscernible sequence $I$,
$R\in\cR$ and a sequence $b$, if $R(b,c)$ holds for at least
$\kappa(M)$ many $c \in I$, then it holds for all but possibly
fewer than $\kappa(M)$ many $c \in I$. This is a property of
stability that we will use in the proof of Lemma~\ref{lem:stab1}.

We will show the following proposition, which is easily obtained
from \cite{sh.hy:strong} and \cite{sh:book} Lemma III 1.11 (3)
given the subsequent lemma. By \emph{Lascar strong types are
stationary in $M$}, we mean that for each finite $a$ and set $A$,
if $b_\ell$ realizes $lstp(a/A)$ and $b_\ell$ is free from $C$
over $A$, for $\ell = 1,2$, then $tp(b_1/C) =tp(b_2/C)$.

\begin{proposition}\label{prop:stab_main}
Let $(M,\cR)$ be simple. The following conditions are equivalent:
\begin{enumerate}[(1)]
\item $M$ is stable. \item Lascar strong types are stationary in
$M$. \item There exists cardinals $\nu \leq \mu <
\beth_{(2^{|\cR|^+})}$ such that for each $A$ and for each $B$
containing $A$, the set
\[
\{\, tp(a/B):\: tp(a/B)\text{ is free over }A\, \}
\]
has size at most $(|A| + \mu)^{<\nu}$.
\end{enumerate}
\end{proposition}

We first show:

\begin{lemma}\label{lem:stab1}
 Suppose that $(M,\cR)$ is simple and stable.
\begin{enumerate}[(1)]
\item If $p$ splits strongly over $A$ then $p$ divides over $A$.
\item If $A$ is free from $C$ over $B$ then $tp(a/B\cup\{c\})$
does not split strongly over $B$ for each finite $a \in A$, and $c
\in C$.
\end{enumerate}
\end{lemma}

\begin{proof}
(1) We show the contrapositive. Suppose that $p$ does not divide
over $A$ and let $I=\{\, a_i:\: i < \mu\,\}$ be indiscernible over
$A$ where $a_0, a_1$ are in the parameters of $p$. Let $R\in\cR$
be such that $R(x,a_0) \in p$. We must show that $R(x,a_1) \in p$.
Consider $p(x,a_0,a_1) = p \restriction Aa_0a_1$. By
$s\lambda-$homogeneity, we may assume that $\mu \geq \kappa(M)$.
Notice that $\{\, a_{2i}\hat{\:} a_{2i+1} :\: i < \mu \,\}$ is
indiscernible over $A$. Then $\bigcup_{i < \mu}
p(x,a_{2i},a_{2i+1})$ is realized by some $b \in M$, since $p$
does not divide over $A$. Hence $R(b,a_{2i})$ holds for all $i <
\mu$, and by stability, $R(b,a_i)$ holds for all but fewer than
$\kappa(M)$ elements $a_i \in I$. This implies that
$R(b,a_{2i+1})$ holds for some $i < \mu$. Hence $R(x,a_{2i}),
R(x,a_{2i+1}) \in p(x,a_{2i},a_{2i+1})$, so in particular
$R(x,a_1) \in p(x,a_0,a_1)$, i.e. $R(x,a_1) \in p$. The proof that
$R(x,a_1) \in p$ implies that $R(x,a_0) \in p$ is similar.

(2) Follows from (1) by definition of freeness (freeness has
finite character).
\end{proof}

\begin{proof}[Proof of the proposition]
Assume that $(M,\cR)$ is simple. (1) follows from (2) or (3) by
counting types just like in the first-order case (we use the fact
that each type is free over a subset of size less than $\kappa$ of
its parameters). To see that (1) implies (2), use the preceding
lemma and Theorem 3.12 of \cite{sh.hy:strong}, which says that if
$tp(a/Ac)$ and $tp(b/Ac)$ do not split strongly over $A$ and have
nonsplitting extensions over any set, and if in addition
$lstp(a/A)=lstp(b/A)$ then $tp(a/Ac) = tp(b/Ac)$. (1) implies (3)
is also as in the first-order case: Let $A$ and $B$ containing $A$
be given. Assume $(M,\cR)$ is stable. We show that (3) holds with
$\nu = \kappa(M)$ and $\mu = \lambda(M)$. By the preceding lemma,
it is enough to show that the number of extensions in $S(B)$ which
do not split strongly over $A$ is at most $\epsilon =(|A| +
\lambda(M))^{<\kappa(M)}$. Let $B' \subset B$ contain a
realization for each $lstp(b/A)$ for $b$ a finite sequence in $B$.
Then, $B'$ has size at most $\epsilon$: By the stability spectrum
theorem, $(M,\cR)$ is stable in $\epsilon$, so any set of size
$\epsilon^+$ contains an infinite $A$-indiscernible subset, and
thus contains different realizations of the same Lascar strong
type over $A$. Notice that if $p$, $q \in S(B)$ do not split
strongly over $A$ and $p \restrict B' = q\restrict B'$, then $p =
q$. To see this, suppose $R(x,b) \in p$. Let $b' \in B'$ realize
$lstp(b/A)$. Then there exist $n < \omega$ and $b_i$ for $i \leq
n$ such that $b_0 = b$, $b'= b_n$ and $b_i, b_{i+1}$ belong to an
infinite $A$-indiscernible sequence. Since $p$ does not split
strongly over $A$, we have $R(x,b_i) \in p$, for $i \leq n$. Hence
$R(x,b') \in p$, so $R(x,b') \in q$, and so $R(x,b) \in q$ by the
same argument applied to $q$. Hence, the number of free extensions
over $B$ is bounded by $|S(B')|$ which is at most $\epsilon$,
since $(M,\cR)$ is stable in $\epsilon$.
\end{proof}

\section{Examples}
\label{sec:examples}

In this section, we give several examples of simple, or simply
stable mathematical structures. Some of these examples clarify the
connection between stability and simplicity of a homogeneous
structure and its first order theory. In particular, stability
does not imply simplicity and the homogeneous models of a stable
first order theory may not be simple. We also give examples which
are simply stable, but whose first order theory is not. Several
examples of groups also illustrate limitations of possible
generalizations; we may not have generics, or large abelian
subgroups, even under strong stability assumptions. Finally, some
of these examples have a natural and immediate description at this
level of generality, but a less amenable first order behavior.

\subsection{Models of simple and stable theories. }

The first place to look for examples of simple
$s\lambda-$homogeneous models is, of course, the homogeneous
models of a simple first-order theory. However, the situation is
not so clear. Saturated models $M$ of simple or stable theories
are clearly simple, or simply-stable homogeneous structures $(M,
\cR)$, where $\cR$ is the collection of definable relations.
Moreover, freeness agrees with nondividing in this case. It is
also clear that a homogeneous model of a stable first-order theory
is stable in the sense of homogeneous structures (as it realizes
fewer types). However, homogeneous models of even a stable theory
may fail to be simple, because they may fail to have the extension
property (see the example after the lemma). We have:








\begin{lemma}
\label{lem:stab_to_simp} Let $M$ be a large uncountable
homogeneous model of a stable theory $T$, $\cR$ the collection of
definable relations on $M$. For $A\cup\{\aa,\bb\}\subset M$,
$tp_{\cR}(\aa/\bb)$ divides over $A$ implies $tp(\aa/\bb)$ divides
over $A$ in $T$. Thus, if $\kappa=\kappa(T)$, $M$ is almost
$\kappa$-simple.
\end{lemma}

\begin{proof}
We are assuming $|M|>|T|^+$, hence $|M|>\kappa$. Suppose
$tp(\aa/\bb)$ does not divide over $A$ in $T$. If
$tp_{\cR}(\bb/A)$ is small, then $tp_{\cR}(\aa/\bb)$ does not
divide over $A$. So, suppose $tp_{\cR}(\bb/A)$ is large and let
$I$ be any indiscernible sequence in $tp_{\cR}(\bb/A)$ of
cardinality $>\kappa$. For $\bb_0\in I$ there is $\aa'$ such that
$q(\xx,\yy)=tp(\aa'\bb_0/A)=tp(\aa\bb/A)$. There is $J\subset I$,
$|J|<\kappa$, such that $\aa'$ is independent from $I'=I\setminus
J$ over $A$. By the stationarity of strong types in a stable
theory, $\aa'$ realizes $\bigcup_{\cc\in I'}q(\xx,\cc)$. Since
indiscernible sequences in a stable theory are indiscernible sets,
$I$ and $I'$ have the same isomorphism type over $A$. Thus, there
is $\aa''\in M$ realizing $\bigcup_{\cc\in I}q(\xx,\cc)$. This
shows that $tp_{\cR}(\aa/\bb)$ does not divide over $A$.
\end{proof}

The following example is due to Shelah. See \cite{HyLe} for
details. This example shows that some homogeneous models of a
stable theory may fail to be simple. Since a homogeneous model of
a stable first order theory is stable as a homogeneous structure,
this example shows also that stability does not imply simplicity
at this level of generality.

The language consists of an infinite number of binary relations
$E_n(x,y)$, for $n < \omega$. The first-order theory $T$ asserts
that each $E_n(x,y)$ is an equivalence relation on the models of
$T$ with an infinite number of equivalence classes, all of which
are infinite. Furthermore, each $E_{n+1}$-class is partitioned
into infinitely many $E_n$-classes. It is easy to see that $T$ is
complete and $\omega$-stable. Let $\bar{M}$ be a large saturated
model of $T$ and $a \in \bar{M}$. Let $M = \bigcup_{n < \omega}
a/E_n$, where $a/E_n$ is the $E_n$-class of $a$. Then $M$ is a
large homogeneous model of $T$ ($\cR$ is the set of first order
definable relations of $M$). $M$ is $\omega$-stable, and so almost
$\omega$-simple. However, for any $b, c \in M$, the type $tp(b/c)$
divides over the empty set. This shows that no type over the empty
set has a free extension, so $M$ fails to be simple.

\subsection{Trees. }

Let $\beta$ be an ordinal and $M=A^{\leq\beta}$ viewed as a tree
structure, where $|A|>|\beta|$. ($A^{\leq\beta}$ is the set of
sequences from $A$ indexed by ordinals $\leq\beta$.) The language
under which the tree is formulated is immaterial, let's say it is
a partial order $\leq$ interpreted by the subsequence relation. It
isn't difficult to show that $M$ is homogeneous. (All elements of
the same length are in the same orbit.) Dividing in $M$ is
understood as follows. Let $c\in M$ be a sequence of length
$\gamma<\beta$ and suppose $c<b$ and $c<a$. If there is a $d$,
$c<d<a$ and $c<d<b$, then $tp(b/a)$ divides over $c$. To see this
consider an indiscernible sequence $I=\{\,d_i:\:i<\omega\,\}$ in
$tp(d/c)$ such that the $d_i$s are distinct. Since the elements of
$I$ must lie on distinct branches in the tree there is no $x\in M$
satisfying $x>d_i\wedge x>d_j$, for $i\not=j$. Thus, the formula
$x>d$ divides over $c$. Continuing for any $a\in M$, $A\subset M$,
let $B\subset A$ be a minimal set such that if $c<a$ and $c\leq
b$, for some $b\in A$, then $c<b$, for some $b\in B$. Then,
$tp(a/A)$ does not divide over $B$.
>From this observation it is easy to show that $M$ is
 $\kappa-$simple for any $\kappa\leq|\beta|$.
$M$ is also stable.

\subsection{A torsion abelian group. }

For $p$ a prime number let $\Bbb{Z}_p$ be the cyclic group with
$p$ elements. Let $\lambda$ be an uncountable cardinal. Let
\[M=\bigoplus_{p\,\mathrm{ prime}}\bigoplus_{\lambda\,\mathrm{
copies}}\Bbb{Z}_p.\] and $\cR$ the class of existentially
definable relations on $M$ in the language with just $+$ and $0$.
For $p$ a prime let $M_p$ be the elements of order $p$ in $M$.
Then, two tuples from $M_p$ with the same quantifier-free type are
in the same orbit of $\mathrm{Aut}(M)$. In other words, the only
invariant structure on $M_p$ is as a vector space over the field
with $p$ elements.
>From this observation it is easy to see that $M$ is homogeneous
and for any large type $r$ over $A$, where $|A|<|M|$, $M$ contains
an infinite $A-$indiscernible sequence of realizations of $r$.
Moreover, for any prime $p$, if $\{a\}\cup A\subset M_p$,
$B\subset A$, and $tp(a/A)$ divides over $B$, then $a$ is linearly
dependent on $A$ over $B$. Given an arbitrary $a\in M$ there are
primes $p_1,\dots,p_n$ and $b_i\in M_{p_i}$ such that
$a=b_1+\dots+b_n$. In this situation, $tp(a/b_i)$ divides over
$\emptyset$. It follows that $M$ is simply superstable. Notice
that the first-order theory of $M$ is only stable.
\medskip

This example illustrates a complication of dealing with simply
stable homogeneous groups, namely, generic types may not exist. A
type $q\in S(A)$ in a simple homogeneous group $G$ is
\emph{generic} if for all $a\in G$, $aq$ is free from $A\cup\{a\}$
over $\emptyset$. For any $q\in S(\emptyset)$ there are primes
$p_1,\dots,p_n$ such that a realization of $q$ is in
$M_{p_1}+\dots+M_{p_n}$. If $b\in M_p$, for $p\not= p_i$,
$i=1,\dots,n$, then $bq$ divides over $b$. Thus, there are no
generic types over $\emptyset$ in $M$.

This context is too general to guarantee the existence of generic
types as in simple first-order theories. In his thesis, A.
Berenstein \cite{be:thesis} considers simply stable homogeneous
groups under the assumption that generics exist and proves that
they behave as in the first order (i.e., compact) case.

\subsection{Free groups.}

Another natural example of simply stable homogeneous groups are
free groups. Let $G$ be an uncountable free group generated by the
subset $S$ of $G$, in the language of groups $(*, ^{-1},1)$, and a
unary predicate $S$ for the set of generators. Then $G$ is a
homogeneous structure (which is not saturated, since the type of
an element not generated by $S$ is omitted). Any two uncountable
such groups of the same size are isomorphic, and it was observed
by H. J. Keisler~\cite{Ke} that these groups are an example of his
categoricity theorem for $L_{\omega_1, \omega}$. Uncountable
categoricity implies that $G$ is $\omega$-stable as a homogeneous
structure (see Section 5); this can also be seen directly by
counting types. $G$ is also supersimple. The freeness relation can
be described easily. For $A \subseteq G$, let $A' \subseteq S$ be
those elements generating the elements of $A$. $A'$ is uniquely
determined from $A$ and this operation commutes with automorphisms
of $G$. For $B \subseteq A,C \subseteq G$, the reader can check
that $\indep{A}{B}{C}$ if and only if $A'$ is free from $C'$ over
$B'$ in the sense of the trivial pregeometry $S$.

This is again a case where we fail to have generics. Moreover,
each abelian subgroup of $G$ is countable, so we cannot expect a
generalization of Cherlin's theorem stating that an uncountable
saturated group whose first order theory is $\omega$-stable has a
definable abelian subgroup of the same size.

\subsection{Hilbert spaces. }

Some of these ideas appear earlier in the work of Henson and
Iovino on Banach space structures, see for example \cite{Io}. The
main technical difference is that the logic of Banach space
structures is set-up so as to keep compactness. They have several
ways of measuring the space of types according to the density
character of these spaces with respect to various topologies.
However, they show that stability is independent of the topology.
Moreover, Iovino showed recently that a compact, homogeneous
Banach space structure is stable in their sense if and only if it
is stable in the sense developed in the previous section. Some of
the material below was also worked out simultaneously, as well as
extended, by A. Berenstein in his Ph.D. thesis~\cite{be:thesis}.
He produced several simply stable expansions of Hilbert spaces and
proved that the $C^*$ algebra of square integrable functions on a
set $X$ is also a simply stable homogeneous structure. He also
observed that, as a group, a Hilbert space does not have generics.
See also \cite{be.bu:hilbert} for more simply stable expansions of
Hilbert spaces.

Let $K$ be the real or complex numbers. An \emph{inner product
space} (or \emph{pre-Hilbert space}) over $K$ is a $K-$vector
space $V$ equipped with a map $\langle-,-\rangle$ from $V\times V$
into $K$ satisfying, for $x,\,y,\,z\in V$ and $r\in K$:
\begin{itemize}
\item $\langle x+y,z\rangle=\langle x,z\rangle+\langle
y,z\rangle$, \item $\langle rx,y\rangle=r\langle x,y\rangle$,
\item $\langle y,x\rangle=\overline{\langle x,y\rangle}$, \item
$\langle x,x\rangle$ is a real number $>0$ if $x\not=0$.
\end{itemize}
If $V$ is an inner product space it is a normed linear space under
 $\parallel -\parallel=\sqrt{\langle x,x\rangle}$.
$V$ is a \emph{Hilbert space} if it is a Banach space under
$\parallel -\parallel$; i.e., $V$ is a complete metric space under
the norm.

For simplicity suppose $K=\Bbb{R}$. This can be formulated as a
model in a first-order language by having one sort, $S_V$, for $V$
and one, $S_F$, for $\Bbb{R}$, a binary map $\langle -,-\rangle$
from $V\times V$ into $\Bbb{R}$, the field operations on
$\Bbb{R}$, $<$ on $\Bbb{R}$, and a function giving scalar
multiplication. For simplicity we assume there is a constant in
the language for each element of $\Bbb{R}$. Let $T_E$ be the
theory in this language expressing the itemized properties above
for the field $S_F$ and that $S_F$ is an ordered field containing
$\Bbb{R}$.

\begin{lemma}
\label{lem:euc_homog} Any pre-Hilbert space $W$ is a subspace of a
Hilbert space $U$ that is homogeneous with respect to
quantifier-free relations.
\end{lemma}

\begin{proof}
Let $M_0$ be a homogeneous universal model of the universal theory
$T_E$ (see \cite{pi:exist} or \cite{sh:lazy} for an earlier
reference) that contains $W$ and has arbitrarily large
cardinality. Let $\cR_0$ be the quantifier-free definable
relations on $M_0$. Let $E(x,y)$ be the equivalence relation on
$M_0$ that holds if $\langle x-y,x-y\rangle< 1/n$, for all
$n<\omega$. Thus, $E$ is $\cR_0-$type definable. Let $(M',\cR')$
be the logical structure on the sort $M_0/E$ in
$(M_0^{eq},\cR_0^{eq})$. By Lemma~\ref{lem:compact_imag},
$(M',\cR')$ is homogeneous. The model $M'$ consists of $V'$
interpreting $S_V$ and $F'$ interpreting $S_F$. The reader should
verify that $M'$ satisfies the axioms for a pre-Hilbert space
relative to the field $F'$ under the operations and relations
inherited from $M$. The real numbers are a subring of $F'$ since
each real is a constant symbol of the language. Notice that if
$s\in F'$, and $|s|\leq r$ for some $r\in\Bbb{R}$, then
$s\in\Bbb{R}$. Let $U=\{\,a\in M':\:\parallel a\parallel
\in\Bbb{R}\,\}$.

\begin{claim}
$U$ is a Hilbert space homogeneous with respect to quantifier-free
relations.
\end{claim}

Let $x,\,y\in U$. In $M'$, $|\langle x,y\rangle|\leq \parallel
x\parallel\cdot\parallel y\parallel\in\Bbb{R}$, hence $\langle
x,y\rangle\in\Bbb{R}$. Thus, $U$ is a submodel of $M'$. The
crucial axioms above are universal, so $U$ is also a pre-Hilbert
space. Since $M_0$ is an uncountable homogeneous universal model
of $T_E$ it realizes any consistent quantifier-free type over a
countable subset. It follows that $\parallel-\parallel$ defines a
complete metric space on $U$. Moreover, it is easy to verify that
$U$ is homogeneous with respect to the quantifier-free relations.
This proves the claim and the lemma.
\end{proof}

Recall the condition $(\bold{P})$ defined at the beginning of
Section~\ref{sec:div}. We assume from hereon that $M$ is a
homogeneous Hilbert space of cardinality $\lambda$, where
$\lambda$ is sufficiently large so that $(\bold{P})$ holds for
$\pi=\aleph_1$ and $\pi'$ some uncountable cardinal $\leq\lambda$.
Here, of course, we mean that $M$ is homogeneous with respect to
quantifier-free relations.

\begin{theorem}
\label{thm:euc_simp} $M$ is $\aleph_1-$simply stable.
\end{theorem}

The proof will be in a series of lemmas. Here are some of the
basic facts about Hilbert spaces that go into the proof. As a
reference we suggest \cite[Ch.16]{he.st:book}.

\begin{remark}\label{rem:hilbert1}
(i) If $H$ is a Hilbert space and $A\subset H$, the minimal
Hilbert space $H'\subset H$ that contains $A$ is denoted
$\hat{A}$. Since $\hat{A}$ is obtained by successively closing
under the vector space operations and taking the limit of Cauchy
sequences, $\hat{A}\subset dcl(A)$; i.e., it is invariant under
the automorphisms of $H$ that fix $A$.

(ii) If $H$ is an inner product space and $A\subset H$ then $A$
contains an orthogonal set $E$ such that every $x\in A$ is
nonorthogonal to some element of $E$ (by Zorn's Lemma). Such an
$E$ is called a \emph{complete orthogonal set for $A$}.
\end{remark}

\begin{notation}
Given an inner product space $H$, $x\in H$ and $Y\subset H$,
$x\perp Y$ if $\langle x,y\rangle = 0$, for all $y\in Y$.
\end{notation}

The following proposition summarizes the elementary results that
enter the proof.

\begin{proposition}
\label{prop:basic_hilbert} (i) (Bessel's Inequality) Let $E$ be a
nonempty orthonormal set in an inner product space $H$, and let
$x\in H$. Then $\sum_{z\in E} |\langle x,z\rangle|^2\leq \parallel
x\parallel^2$, hence $\{\,z\in E:\:\langle x,z\rangle\not=0\,\}$
is countable.

(ii) Let $H$ be a Hilbert space and $E$ a complete orthogonal
subset of $H$. Then for all $x\in H$, $x=\sum_{z\in
E}\frac{\langle x,z\rangle}{\langle z,z\rangle}z$.
\end{proposition}

\begin{lemma}
\label{lem:ortho_dnd} Given $A\cup E\subset M$, with $E$ an
orthogonal set over $A$, and $a\in M$ such that $a\perp E$,
$tp(a/A\cup E)$ does not divide over $A$. If $B\subset dcl(A\cup
E)$, then also $tp(a/A\cup E\cup B)$ does not divide over $A$.
\end{lemma}

\begin{proof}
Let $\{\,E_i:\:i\in X\,\}$ be an $A-$indiscernible sequence in
$tp(E/A)$. The construction of $M$ at the beginning of this
section, and the fact that $E$ is orthogonal over $A$ guarantees
the existence of an $a'$ realizing $tp(a/A)$ such that $a'\perp
E_i$, for all $i\in X$. Checking the possible relations on $A\cup
E\cup \{a\}$, this implies that
$tp(E_i\cup\{a'\}/A)=tp(E\cup\{a\}/A)$, for all $i\in X$. This
witnesses that $tp(a/A\cup E)$ does not divide over $A$.

Continuing, for each $i\in X$, there is $B_i$ such that
$tp(E_i\cup B_i/A)=tp(E\cup B/A)$. Since $B_i\subset dcl(E_i\cup
A)$, $tp(B_i/E_i\cup A)$ has a unique extension over $A\cup
E_i\cup\{a'\}$. Thus, for each $i\in X$, $tp(B_i\cup
E_i/A\cup\{a'\})=tp((B\cup E/A\cup\{a'\})$. We conclude that
$tp(a/A\cup E\cup B)$ does not divide over $A$.
\end{proof}

\begin{lemma}
\label{lem:free1} Suppose $a\in M$ and $A\subset M$, where
$|A|<|M|$. Then there is a countable $B\subset A$ such that
$tp(a/A)$ does not divide over $B$.
\end{lemma}

\begin{proof}
Let $E$ be a complete orthogonal subset of $\hat{A}$, chosen so
that $E\cap A$ is a complete orthogonal subset of $A$. Let
$F=\{\,z\in E:\:\langle a,z\rangle\not=0\,\}$, which is countable
by Bessel's Inequality. Let $F_0=F\cap A$ and $F_1=F\setminus
F_0$. Let $A_0$ be a minimal subset of $A$ such that $F_1\subset
dcl(A_0)$. Since $F$ is countable, $A_0$ is countable by
Remark~\ref{rem:hilbert1}(i), and we may assume that $F_0\subset
A_0$. By enlarging $A_0$ and $F$ if necessary we can assume that
$e\in E\setminus F$ implies that $e\perp A_0$. By
Lemma~\ref{lem:ortho_dnd}, $tp(a/\hat{A})$ does not divide over
$A_0$, proving the lemma.
\end{proof}

\begin{lemma}
\label{lem:free_ext} Given $a\in M$ and $A\subset B\subset M$
there is a unique complete type over $B$ which is an
$\aleph_0-$free extension of $tp(a/A)$, and this type does not
divide over $A$.
\end{lemma}

\begin{proof}
Consider the Hilbert spaces $\hat{A}\subset \hat{B}$. Let $E$ be a
complete orthogonal subset of $\hat{A}$ and $F\supset E$ a
complete orthogonal subset of $\hat{B}$. Notice that each $f\in
F\setminus E$ is orthogonal to $\hat{A}$. Since $\hat{A}\subset
dcl(A)$, $tp(a/A)$ has a unique extension over $\hat{A}$, namely
$tp(a/\hat{A})$, which also does not divide over $A$. There is a
$b$ realizing $t(a/\hat{A})$ such that each $f\in F\setminus E$ is
orthogonal to $b$, and among such $b$ there is a unique type over
$\hat{A}\cup F$. There is a unique extension of $tp(b/\hat{A}\cup
F)$ over $\hat{B}$, and by Lemma~\ref{lem:ortho_dnd}
$tp(b/\hat{B})$ does not divide over $A$. This proves the lemma.
\end{proof}

Combining Lemma~\ref{lem:free1} and Lemma~\ref{lem:free_ext} shows
that $M$ is $\aleph_1$-simple and stable, proving
Theorem~\ref{thm:euc_simp}.

\bibliography{papersbib}
\bibliographystyle{alpha}

\end{document}